\documentclass[12pt,english]{article}

\usepackage{graphicx,psfrag,epsfig,amsfonts,amssymb,amsmath,babel}

\newcommand{\dist}{\mbox{$\,$dist$\,$}}

\newtheorem{lemma}{Lemma}[section]
\newtheorem{theorem}[lemma]{Theorem}
\newtheorem{proposition}[lemma]{Proposition}

\newtheorem{teorema}[lemma]{Theorem}

\newtheorem{corollary}[lemma]{Corollary}

\newtheorem{definition}[lemma]{Definition}

\newtheorem{remark}[lemma]{Remark}

\newtheorem{example}[lemma]{Example}

\newtheorem{prueba}[lemma]{}
\begin{document}

\title{Observable invariant measures.}
\author{Eleonora Catsigeras and Heber Enrich
 \thanks{
  Both authors: Instituto de Matem\'{a}tica y Estad\'{\i}stica Rafael Laguardia (IMERL),
 Fac. Ingenier\'{\i}a,  Universidad de la Rep\'{u}blica,  Uruguay.
 E-mails: eleonora@fing.edu.uy
 and enrich@fing.edu.uy Address: Herrera
  y Reissig 565. Montevideo. Uruguay. }}

\date{February 18th., 2010.}
\maketitle

\begin{abstract}

  For  continuous maps on a compact manifold $M$, particularly for
 those  that do not  preserve the
Lebesgue measure $m$, we define the \em observable invariant
probability measures \em as a generalization of the physical
measures. We prove that   any continuous map has observable measures,
and characterize those that are  physical  in terms of the
observability.  We prove that there exist   physical measures whose basins cover
 Lebesgue a.e, if and only if the set of all observable measures is finite or infinite numerable.
  We   define  for any continuous map,
its \em generalized attractors \em using the set of observable
invariant measures where there is no physical measure, and prove
that any continuous map  defines a decomposition of the space in up
to infinitely many generalized attractors whose basins cover
Lebesgue a.e. We apply the results to the $C^1$ expanding maps $f$
in the circle, proving that the set of all observable measures (even
if $f$ is not $C^{1 + \alpha}$) is a subset of the set of all the
equilibrium states of $-\log |f'|$.

\end{abstract}

\section{Introduction}

It is an old problem to find  \lq\lq good"\ probability measures for a map
$f\colon M \mapsto M$, meaning for that, an invariant probability that resume in some
sense, the dynamics by  iterations of the map. Sometimes, the map is born
with a good measure, as in the case of  maps preserving a Lebesgue ergodic measure. 
But this is not
true in general, and it is not an easy question to determine, in
most examples, a single or a few probability measures representing the dynamics of
the map.

There have been proposed several ideas to define a \lq\lq good"\
invariant probability measure~$\mu$:
\begin{enumerate}
\item \label{phys} Lebesgue a.e. point in a set is generic with respect
$\mu$, that is, $\displaystyle \mu =\lim_{n\to
\infty}\frac1n\sum_{j=0}^{n-1} \delta _ f^j(x)$, where the
convergence is in the weak$^*$ topology of the space ${\cal P}$ of
probabilities on $M$.
\item \label{2srb} The conditional measures of $\mu$ on unstable manifolds
are absolutely continuous respect to the Lebesgue measure along those manifolds.
\item \label{pesin} $\mu$ verifies the Pesin-Ledrappier-Young (PLY) equality:
$$h_\mu (f) =\int \sum_i
\lambda^+_i(x) \dim E_i(x) \, d\mu(x) $$ where $h_{\mu}(f)$ is the entropy of $\mu$,
and $\dim E_i(x)$ is the
multiplicity of the positive Lyapunov exponent $\lambda^+(x)$ in the Oseledec's decomposition.
\item \label{stoch} The measure is the limit of measures which are invariant
under stochastic perturbations.
\end{enumerate}

It is a remarkable property that the four definitions above are equivalent for Axiom A
attractors. Moreover,  Ledrappier and Young
(see \cite{ledrappieryoung} and \cite{ledrappieryoung2}), under
suitable hypothesis of differentiability, proved that a measure
verifies property~\ref{2srb} if it verifies property \ref{pesin}
while the converse result is the well known result of Pesin's entropy formula.
Ergodic measures verifying~\ref{2srb} (or, equivalently,
\ref{pesin}) with no zero Lyapunov exponents describe chaotic
behavior, and are accompanied by rich geometric and dynamical
structures
  (see \cite{young}). Nevertheless, the
other measures listed above are also interesting, because reveal
statistical aspects of the  behavior of the future iterations of the
map. For instance,   a measure verifying  definitions \ref{phys} or \ref{stoch}
is concentred in the part of the space which is statistically
more visited.

 We will call physical measure, a probability measure
verifying \ref{phys}, and stochastically stable, a
probability measure verifying \ref{stoch}. We will call SRB
(Sinai-Ruelle-Bowen) measure a probability  verifying \ref{2srb}, and
 a PLY  measure, a probability
verifying~\ref{pesin}.
In this work, we   propose another
 concept of \lq\lq good" probabilities, which we call   observable measures. The
following question  was the   motivation of this work: Is it
possible to describe probabilistically in the space, in some \em
minimal way, and in a very general regular or irregular setting,
\em the asymptotic    behavior of the time averages of
Lebesgue a.e. orbit? We      answer   this  question in Theorem~\ref{toeremaminimal0}.

Generalized ergodic attractors and observable measures, that we
define and theoretically develop along this work, do always exist
for any continuous map (Theorems \ref{teoremaexistencia0} and
\ref{teoremaDescomposicionGenErgAttr0}). On the other hand, physical
measures and ergodic attractors  do not necessarily exist (see examples
\ref{ejemplodossillas} and \ref{ejemploExpanding}).

It is largely known the difficulties to
characterize, or just find,  non hyperbolic or non $C^{1+ \alpha}$
maps that do have physical measures. This is a hard problem even in
some systems whose iterated topological behavior is  known
(\cite{carvalho}, \cite{enrich}, \cite{hu}, \cite{huyoung}). The
difficulties appear  when trying to apply to a non hyperbolic
setting, or to a non $C^{1+ \alpha}$ map,
 the known  techniques for
constructing the physical measures of hyperbolic $C^{1 + \alpha}$ maps (\cite
{pesin}, \cite {sinai}, \cite {anosov}). The $C^{1 + \alpha}$
hypothesis allows the existence of SRB measures (\cite{bowen3},
\cite{ruelle}, \cite{sinai}) and are relevant and widely studied
 occupying an important focus of
interest in the ergodic differentiable theory of dynamical systems
(\cite{anosov},  \cite{pesinsinai}, \cite{pughshub2},
\cite{viana}, \cite{bonattiviana}).
But if
having a weak or non uniform hyperbolic setting,  the obstruction usually resides in the
irregularity of the invariant manifolds, which technically
translate the non trivial relations between topologic, measurable
and differentiable properties of the system (\cite{pughshub2}).
The difficulties arise  even
for maps posed in a very regular setting as Lewowicz diffeomorphisms in the
two-torus (see \cite{lewowicz}).  For them, the differentiable
regularity of the given transformation (they are analytic maps)
and the topological  known behavior of the  iterated system (they
are conjugated to Anosov maps), were not enough, up to the moment, to prove the
existence of physical measures and ergodic attractors, except in some meager set of
examples \cite{nosotros}.
Before stating the results we need to formalize some definitions
that we will use all along this work: Let $f\colon M \mapsto M$ be a
continuous map with $M$ a compact, finite-dimensional manifold. Let
$m$ be a probability Lebesgue measure, and not necessarily
$f$-invariant. We denote $\cal P$ the set of all Borel probability
measures in $M$, provided with the weak$^*$ topology, and a metric
structure inducing this topology.

 For any point $x \in M $ we denote
$p\omega (x)$ to the set of the Borel probabilities in $M $ that are
the partial limits of the (not necessarily convergent) sequence
\begin{equation} \label{1} \left\{ \frac{1}{n} \; \sum_{j=0}^{n-1} \delta
_{f^j(x)} \right\} _{ n \in \mathbb{N}} \end{equation} where
$\delta_y$ is the Dirac delta probability measure supported in $y
\in M$.

 The set
$p\omega (x) \subset \cal P$ is the collection of the spatial
probability measures describing the asymptotic time average (given
by (\ref{1})) of the system states, provided the initial state is
$x$. If the sequence (\ref{1}) converges then we denote $p\omega (x)
= \{\mu _x\}$. To include also those cases for which
the sequence (\ref{1}) is not convergent (for a set of orbits with
positive Lebesgue measure) we consider, for a given measure $\mu$,
the set of points $x \in M$ such that the minimum distance between
$\mu$ and the set of partial limits of the sequence (\ref{1})  is small. We define:

\begin{definition}
 \em \label{definicionobservable} \label{definicion1} {\bf (Observable probability measures.)}
A probability measure $\mu \in \cal P$ is \em observable \em if for
all  $\varepsilon
>0 $  the set ${A}_{\varepsilon }= \{x \in M:
 \mbox{ dist}^* ( p\omega(x), \mu )<  \varepsilon \}$ has
 positive Lebesgue measure.
The set ${A}_{\varepsilon}  = A_{\varepsilon}(\mu)   \subset M$ is
called \em the $\varepsilon$- basin of partial attraction  \em of
the probability $\mu$.
\end{definition}

We note that the definition above is independent of the choice of
the distance  in $\cal P$, provided that  the  metric structure
induces its weak$^*$ topology. Observable
measures  are $f$-invariant, and  usually at most a few part
of the space  of invariant measures for $f$ are observable
measures (see the examples  in Section \ref{ejemplos2}).
%In paragraph \ref{interpretacionfisica} of the Appendix, we explain
%the mathematical and physical reasons that lead us to use the
%$\varepsilon >0$ condition in Definition \ref{definicionobservable}
%instead of using $\varepsilon = 0$ at once.
 We remark that for observable measures,
 the condition $m(A_\varepsilon) >0$ must be verified not only
 \em for some \em  but \em for all \em $\varepsilon >0$.

\begin{definition}
 \em \label{definicionsrb} {\bf (Physical probability measures.)}
A probability measure $\mu \in \cal P$ is \em physical \em (even if it is not ergodic), if
 the set $B = \{x \in M:
  p\omega(x)= \{\mu\}  \}$ has
 positive Lebesgue measure.
The set  $B= {B} (\mu) \subset M$
  is called the \em basin of attraction \em of $\mu$.
\end{definition}

Therefore, all physical measures are
observable, but not all observable measures are physical (Examples \ref{ejemplodossillas} and
\ref{ejemploExpanding}).

  \vspace{.3cm}

We state  the following starting results:

\begin{teorema} \label{teoremaexistencia0} \em  { \bf (Existence of
observable measures and physical measures.) } \em

{\bf a)} For any continuous map $f$, the set $\cal O$ of all
observable probability measures for $f$ is non-empty and weak$^*$-
compact.

{\bf b)} $\cal O$   is finite or countably infinite if and only if
there exist \em (resp. finitely or countable infinitely many) \em
physical measures of $f$ such that the union of their basins of
attraction cover  Lebesgue a.e.
\end{teorema}

The first statement  of this theorem is proved in paragraph
\ref{teoremaexistencia} and the second one in paragraph \ref{srb}.

\vspace{.3cm}

The $p\omega(x)$ limit set of convergent subsequences of (\ref{1})
may have many different partial limit measures. Nevertheless, we
prove that
 $p\omega (x)$ is formed only
with observable measures, for Lebesgue almost all $x \in M$, as stated in
the following Theorem \ref{toeremaminimal0}.

\begin{definition} \em \label{definicionbasin}  {\bf
(Basin of attraction.)}

 The \em basin of attraction $B({\cal K})$ \em  of a compact
subset ${\cal K}$ of the space ${\cal P}$ of all the Borel
probability measures in $M$, is the (maybe empty) subset of $M$
defined as:
$$B({\cal K}) = \{x \in M: p\omega(x) \subset {\cal K} \}$$
\end{definition}

If the purpose is to study the asymptotic to the future time
average behaviors of Lebesgue almost all points in $M$, then the
\em consideration of the set $\cal O$ of all the observable measures for $f$,
 is the necessary  and
sufficient condition. \em   In fact we have the following:
\begin{teorema} \label{toeremaminimal0} \em  {\bf
(Attracting minimality property of the set of observable measures.)} \em
The set $\cal O$ of all observable measures for $f$ is the minimal
 compact subset of the space $\cal P$ whose basin of attraction has
total Lebesgue measure.
\end{teorema}
We prove this theorem in paragraph \ref{pruebaGenErgAttract}.

\vspace{.3cm}

The theory about the observable measures
 describe  the statistical asymptotic behavior of
time averages,  instead   of the theory of physical measures  when
these last probabilities do not exist. In particular it is suitable
to study the statistics of  the future iterations of maps, disregarding their regularity,
 that do not preserve a Lebesgue equivalent measure,  or that do preserve it but
 are not ergodic. In
other words,  the results about the observable measures,
independently if they are or not useful to find physical measures in
some concrete examples, substitute the   physical measures, and do that in
a wide setting of dynamical systems (all the continuous maps),
without loosing their   statistical meaning (see Proposition
\ref{propositionAtraeEnMedia}.)

%In resume, the    physical
%or natural measures are widely studied
% occupying a relevant
%interest in the ergodic theory of dynamical systems,
%but very few results are known if the system is not more than $C^{1 + \alpha}$,  and if the
%technique of bounded distortion to find absolute continuous probabilities fails.
%

%In this paper we characterize the physical measures as a particular
%class of observable measures. In fact, if they exist, all physical measure are observable, and conversely, if there exists some observable measure that is isolated in the space of all observable measure, then it is physical (see Theorem AGREGAR REFERENCIA).

\vspace{.3cm}

 In Section \ref{seccion2} we define
\em  attractor $A$, \em  to the support in $M$ of a physical, non
necessarily ergodic measure $\nu$, similarly as done by  Pugh and
Shub in \cite{pughshub}. Analogously, we call  \em generalized
attractor, \em  to the union $A$ of the  supports of an adequately
reduced weak$^*$-compact family ${\cal O}_1$  of observable
measures. We will not require strict topological attraction to $A$, but
weak topological: in the average   the orbit  lays in an arbitrarily
small neighborhood of $A$ as much time as wanted. (Proposition
\ref{propositionAtraeEnMedia}.) In Theorem
\ref{teoremaDescomposicionGenErgAttr0} we prove that there always
exist up to countable many generalized attractors whose basins cover
Lebesgue a.e. The following open question refers to the existence
and finiteness of physical measures and to the convergence of the
sequence (\ref{1}) of time averages of the system for a set of
initial states with total Lebesgue measure. It is possed in
\cite{palis}:

\begin{prueba} \em \label{Palis1}
{\bf Palis Conjecture }  \em Most dynamical systems have up to
finitely many physical measures (or ergodic attractors) such  that
their basins of attraction cover Lebesgue almost all points. \em
\end{prueba}

This conjecture admits the following equivalent statement, that
seems weaker. (In fact,  the definition \ref{definicionobservable}
of observability  is certainly weaker than the definition
\ref{definicionsrb}
 of physical measures.)

\begin{prueba} \em \label{Palis2}
{\bf Equivalent formulation of Palis Conjecture:}

\em For most
 dynamical systems the set of observable measures is finite. \em
\end{prueba}

Note: To prove the equivalence of statements \ref{Palis1}  and
\ref{Palis2} it is enough to join the results of  Theorems \ref{teoremaexistencia0}.b
and \ref{toeremaminimal0}.

\vspace{.2cm}

 For systems preserving the Lebesgue measure the main question is
 their ergodicity. It is immediate
   the following result:

\begin{remark}
\label{teoremaLebesgue0} \em  {\bf (Observability and
ergodicity.)}

If $f$ preserves the Lebesgue measure $m$ then the following
assertions are equivalent:

\em (i) \em  $f$ is ergodic respect to $m$.

\em (ii) \em There exists a unique observable measure $\mu$ for
$f$.

\em (iii) \em  There exists a unique physical measure $\nu $ for $f$
attracting Lebesgue a.e.

Besides, if  the assertions above are verified, then $m = \mu =
\nu$

\end{remark}

 Given a  Lebesgue measure preserving map $f$,
 the
 question if $f$ has  a physical measure is mostly open,  for differentiable
 maps that do not have some kind of uniform total or partial hyperbolicity \cite{viana}.
 The existence
 of a  physical measure attracting
 Lebesgue a.e. is equivalent to the ergodicity of the map, and is also an open
question if most of conservative maps are ergodic (\cite{pughshub2}, \cite{bonattivianawilk}). The
key difficult point  resides in those maps that do not have any kind of
uniform total or partial $C^1$ hyperbolicity in the space, as in the  inspiring Lewowicz examples in the two torus
\cite{lewowicz}.  Due to Remark
\ref{teoremaLebesgue0} those open questions   are
  equivalent to the \em uniqueness of the   observable measure. \em
  %We do not know if it helps to solve the question, but at least, note that the definition
%  of observability of a  measure \em is   weaker \em than the definition of its
%  ergodicity, and it is \em always verified \em for some (invariant) measure.

\section{Attractors} \label{seccion2}

 Due to the conjecture  in \ref{Palis1} and Theorem \ref{teoremaexistencia0},
  we are
interested in partitioning the set $\cal O$ of observable
measures, or to reduce it as much as possible, into different
compact subsets whose basins of attractions have positive Lebesgue
measure. Due to Theorem \ref{toeremaminimal0}, no
proper compact part of $\cal O$ has a total Lebesgue
basin. We define:

\begin{definition} \label{definicionreduccion0} \em {\bf
(Generalized  Attractors - Reductions of the space ${\cal O}$.)} A
\em generalized attractor $(A, \cal A) \subset M \times {\cal O}$,
\em (or a \em reduction \em ${\cal A}$ of the space ${\cal O } $ of
all observable measures for $f$), is a compact subset $(A, \cal A)$
such that the basin of attraction $B({\cal A})= \{x \in M: p \omega
(x) \subset {\cal A}\}$  has positive Lebesgue measure in $M$, and
$A $ is the (minimal) compact support in $M$ of all the probability
measures in ${\cal A}$. We call $(A, \{\mu\})$ an \em attractor \em
if it is a generalized attractor with a single invariant probability
$\mu$, i.e. ${\mu}$ is a physical measure.
\end{definition}

\begin{remark}
\em Sometimes we refer to a generalized attractor only to $A$ or only to ${\cal A}$.
The \em irreducible \em generalized attractors (if they exist), attract
  the time averages distributions  and are minimal in some sense, but
 are not formed necessarily
with ergodic measures (Example
\ref{ejemplodossillas}). In spite a system could not exhibit a physical measure, still the
reductions of the space of observable measures divide the manifold
in the basins of  \em generalized attractors. \em Each
reduction ${\cal A}$ has a basin $C = C({\cal A})$ with positive
Lebesgue measure and is minimal respect to $C$. We state this
result as follows:

\end{remark}

\begin{teorema} \em  {\bf
(Minimality of generalized attractors.)} \em
\label{teoremaminimal0r}
Any generalized attractor ${\cal A}$
is the minimal compact set of observable measures attracting its
basin $C({\cal A})$. More precisely
$$m (C({\cal A}) \setminus C({\cal K}) ) >0 \ \ \ \mbox{
for all compact subset } \ \ \ {\cal K} \subset {\cal A}.$$
\end{teorema}

We prove Theorem  \ref{teoremaminimal0r}  in paragraph
\ref{pruebaGenErgAttract}.

\begin{definition}
 \em {\bf Irreducibility } \label{definicionirreducible0}

 A generalized attractor ${\cal A}\subset {\cal P} $ is  \em irreducible \em if it does not
 contain proper compact subsets that are also generalized  attractors.

 It is \em trivial \em or \em trivially irreducible \em
 if its diameter in ${\cal P }$ is zero , or in
  other words, if ${\cal A}$ has a unique observable measure
  $\mu$.

  In other words: physical measures are trivially irreducible and conversely.

\end{definition}

\vspace{.2cm}

 The following
result is  much weaker than, but related with, the  Palis' conjecture
stated in paragraph \ref{Palis1}:

\begin{teorema} \label{teoremaDescomposicionGenErgAttr0} {\bf (Decomposition
Theorem)}

For any continuous map $f\colon M \mapsto M$ there exist a
collection of (up to countable infinitely many) generalized
attractors whose basins of attraction are pairwise Lebesgue-almost
disjoint and cover Lebesgue-almost  all $M$.

The space of all continuous maps divide in two disjoint classes:

$\bullet$ The  generalized attractors of the decomposition are all
irreducible and then the decomposition is unique.

$\bullet $ For all $\varepsilon >0$ there exists   a decomposition
for which the  reducible generalized attractors have all diameter
(in the weak$^*$ space of probabilities), smaller than $\varepsilon$
\em
\end{teorema}
We prove  this theorem in paragraph
\ref{teoremadescomposicionGenArgAttract}.
Note that, in the second class of systems, as the reducible generalized attractors have
in the space of probabities a small diameter, for a rough observer, each of those
attractors    acts as a physical measure.
\vspace{.2cm}

 In Theorem \ref{teoremaPhysicalMeasuresAndChains} and Corollary
\ref{corolarioCo-cadenaSRB} we characterize those maps whose
generalized attractors  are
the support of physical measures, as asked in the statement of Palis
conjecture.

 %Generalized ergodic

%\begin{prueba} \em {\bf A note about the generalized   Attractors and the Milnor
%Attractors.} \label{Milnor2}

\begin{definition} \em {\bf Milnor attractor.} \cite{milnor} \label{definicionMilnor}
\label{Milnor2}
 A Milnor attractor is a compact set $A \subset M$ such that its
 topological basin of attraction $$B(A) = \{ x \in M: \omega (x) \subset A \}$$ has positive
 Lebesgue measure, and for any compact proper subset $K \subset A $
 the set $$B(A) \setminus B(K) = \{x \in M: \omega (x) \subset A ,\;  \omega (x)
 \not \subset K \}$$ also has positive Lebesgue measure.

Note: Here $\omega (x)$ is the $\omega$-limit set in $M$ of the
orbit of $x$, i.e. the set of limit points in $M$ of the orbit
with initial state $x$.
\end{definition}

  The generalized attractors $(A, {\cal A})$ where ${\cal A} \subset {\cal P}$,
    were inspired in the definition
\ref{definicionMilnor}. They play in ${\cal P}$ a similar
topological role that Milnor attractors play in $M$. In particular,
 the physical measures in $\cal P $ play the role
  that sinks do in $M$, those first considered as punctual attractors
  of time averages probabilities, and these last considered as
punctual attractors of the points in $M$, along the future orbits.

\section{Proofs of Theorems \ref{teoremaexistencia0} (a) and \ref{toeremaminimal0}.}

\begin{definition} \em {\bf (Weak$^*$ topology in the space of probability measures.) }
\label{definicion0debilestrella}

 The weak$^*$ topologic structure
in the space $\cal P$ is defined as: $$\mu_n \rightarrow \mu \; \;
\mbox{in } {\cal P} \;\; \; \mbox{ iff } \;\; \; \lim \int \phi \, d
\mu _n = \int \phi \, d \mu \;\; \mbox{ for all } \phi \in C^0(M, \mathbb{R})
$$ where $C^0(M, \mathbb{R})$ denotes the space of the continuous real functions
in $M$.
\end{definition}
A classic basic theorem on Topology states that the space $\cal P$
is compact and metrizable when endowed with the weak$^*$ topology
\cite{mane}.
Let us  denote as ${\cal P}_f \subset {\cal P}$  the set of the
Borel probability measures in $M$ that  are $f$-invariant, that is
$\mu ({f^{-1}(B)}) = \mu (B)$ for all Borel set $B \subset M$. Note
that the Lebesgue measure $m$ does not necessarily belong to ${\cal
P}_f$. Fix any metric in  $\cal P$ giving its weak$^*$ topology
structure. We denote as ${\cal B}_\varepsilon (\mu)$ to the open
ball in $\cal P$ centered in $\mu \in \cal P$ and with radius
$\varepsilon
>0$.

\begin{prueba} \label{definicion0} \em {\bf The $p \omega$- limit sets.}

 At the beginning of this paper we defined, for each initial state $x \in M$,
 the
 set $p\omega (x)$ in the space $\cal P $  as the partial limits, in the
weak$^*$ topology in $\cal P$, of the
sequence (\ref{1}) of time averages.
In other words:
$$p\omega (x) = \left\{\mu \in {\cal P}: \lim \frac{1}{n_i} \sum_{j=0}^{n_i-1}
\phi({f^j(x)}) = \int \phi \, d \mu \; \forall \;  \phi \in C^0(M,
\mathbb{R}) \; \mbox { for some } n_i \rightarrow + \infty\right\}$$
For further uses we state here  the following property for the
$p\omega$-limit sets:
\end{prueba}
\begin{theorem} \em \label{teoremaconvexo} {\bf Convex-like
property.}

 For any point $x \in M$

i) \em If $\mu, \nu \in p\omega(x)$  then for each real number $ 0
\leq \lambda \leq 1$ there exists a measure $\mu_{\lambda} \in
p\omega (x)$ such that $$\dist(\mu_{\lambda}, \mu) = \lambda
\dist(\mu, \nu)$$

 ii) \em The
set $p\omega(x)$ \em either has a single element or non-countable
infinitely many.
\end{theorem}

{\em Proof:} See \ref{prueba} in the appendix.

%\begin{definition}  \label{definicion1} \em
%Given a real number $\varepsilon >0$ we say that a (non-necessarily
%$f$-invariant) probability measure   $\mu$ in $M$  \em is
%$\varepsilon$-observable \em if the set of points $x \in M$ such that
%$p\omega (x) \cap B_{\varepsilon }(\mu) \neq \emptyset$ has positive
%Lebesgue measure.
%\end{definition}

\begin{prueba} \em \label{teoremaexistencia}
{\bf Proof of Theorem \ref{teoremaexistencia0} (a): } (The proof of
the assertion (b) of Theorem \ref{teoremaexistencia0}  is delayed
until paragraph  \ref{srb}.)

Let us prove that, given any continuous map $f$, the set $ {\cal O
}_f $ of the probability measures that are $\varepsilon$-observable
for all $\varepsilon
>0$
 is a
non-empty compact subset of  the (weak$^*$ topologic) space ${\cal
P}_f$ of the $f$-invariant measures.

The key question is that ${\cal O}_f$ is not empty, which we prove
at the end.

 Let us first prove that ${\cal O}_f \subset {\cal P}
_f$. In fact, given $\mu \in {\cal O}_f$ then, for any $\varepsilon
= 1/n
>0$ there exists some $\mu _n \in {\cal B}_{\varepsilon} (\mu)$ which is the
limit of a convergent subsequence of (\ref{1}) for some $x \in M$.
As the limits of all convergent subsequences of (\ref{1}) are
$f$-invariant, we have that $\mu _n \in {\cal P}_f \subset {\cal P}$
for all natural number $n$, and $\mu _n \rightarrow \mu$ with the
weak$^*$ topologic structure of $\cal P$. The space ${\cal P}_f$ is
a compact subspace of $\cal P$ with the weak$^*$ topologic
structure, so $\mu \in {\cal P}_f$ as wanted.

Second, let us prove that ${\cal O} = {\cal O }_f$ is compact. The
complement ${\cal O}^c$ of $\cal O$ in $\cal P$ is the set of all
probability measures $\mu$(not necessarily $f$-invariant) such that
for some $\varepsilon = \varepsilon (\mu) >0$ the set $\{x \in M:
p\omega(x) \cap {\cal B}_{\varepsilon}(\mu) \neq \emptyset\}$ has
zero Lebesgue measure. Therefore ${\cal O}^c$ is open  in $\cal P$,
and ${\cal O}$ is a closed subspace of  $\cal P$. As $\cal P$ is
compact we deduce that  $\cal O$ is compact as wanted.

Third, let us prove that $\cal O$ is not empty. Suppose by
contradiction that it is empty. Then ${\cal O}^c = {\cal P}$, and
for every $\mu \in {\cal P}$ there exists some $\varepsilon =
\varepsilon (\mu) >0$ such that the set $A= \{x \in M: p\omega(x)
\subset ({\cal B}_{\varepsilon}(\mu))^c \}$ has total  Lebesgue
probability.

As $\cal P$ is compact, let us consider a finite covering of $\cal
P$ with such open balls ${\cal B}_{\varepsilon}(\mu)$, say $B_1,
B_2, \ldots B_k$, and their respective sets $A_1, A_2, \ldots A_k$
defined as above. As $m(A_i) = 1$ for all $i= 1, 2, \ldots, k$ we
have that the intersection $B= \cap_{i=1}^k A_i$ is not empty. By
construction, for all $x \in B$ the $p\omega$-limit of $x$ is
contained in the complement of  $B_i$ for all $i = 1, 2 \ldots, k$,
and so it would not be contained in $\cal P$, that is the
contradiction ending the proof. $ \Box $
\end{prueba}

Note that the proof of Theorem \ref{teoremaexistencia0}  (a)
does not use any property of the Lebesgue measure on $M$ different
from those that any Borel probability $m$ on $M$ also has. The same
result works (but maybe defining a different subset of observable
measures) if any other, non necessarily invariant probability measure
$m$ in ${\cal P}$ is used, as the reference probability distribution
for the choice of the initial state $x$, instead of the Lebesgue
measure. If so, the concept of \em physical \em measure also changes
accordingly. Nevertheless, along this work, we are calling $m$ to
the Lebesgue measure, i.e. the volume form measure, given by the
Riemannian metric on the manifold $M$, adequately rescaled to be a
probability: $m(M) = 1$.

\begin{definition}\em \label{definicionObservSize} {\bf Observability size.}
If $\mu \in {\cal P}$ is a (non necessarily invariant) probability
measure  (in particular if it is an observable measure, see
Definition \ref{definicion1}), we call   \em observability size of
$\mu$ \em to the non negative real function $o = o_{\mu}\colon
\overline{\mathbb{R}^+} \mapsto \overline{\mathbb{R}^+}$ defined as
$$o_{\mu } (\varepsilon)= m (A(\varepsilon ,\mu ))$$ where $m$ is the
Lebesgue measure in $M$ and $A(\varepsilon ,\mu )$ is the set
$$A(\varepsilon ,\mu ) = \{x \in M: p\omega (x) \cap {\cal B}_{\varepsilon }(\mu )
\neq \emptyset \}   \; \; \; \; \; \ \mbox{ being } \; \; \; \ {\cal
B}_{\varepsilon }(\mu) = \{\nu \in {\cal P}: \mbox{ dist}^*(\nu,
\mu) < \varepsilon\}$$

For some fixed $\varepsilon >0$, we say that $\mu \in {\cal P}$ is
$\varepsilon-$ observable if $o_{\mu }(\varepsilon) >0$.
\end{definition}

\begin{remark} \em
For any probability measure $\mu$, its observability size function
$o (\varepsilon )$ \em is  positive and decreasing with $\varepsilon
>0$. \em Then  $o(\varepsilon )$  has always a   non-negative limit  value
when $\varepsilon \rightarrow 0^+$.
We reformulate Definition \ref{definicionobservable} of
observability of measures, in the following equivalent terms:
\end{remark}
\begin{prueba} \label{definicion2} \em {\bf Remark: }   {\bf
(Observability revisited.)}
 $\mu \in {\cal P} $ is \em observable \em for $f$, if
and only if it is $\varepsilon$-observable (see \ref{definicion1})
for all $\varepsilon
>0$. In particular, $\mu$ is physical if and only
if $\lim _{\varepsilon \rightarrow 0} o_{\mu}(\varepsilon) >0$.

 The characterization of those continuous
maps having physical-measures as those whose  sets of observable
measures, or some reductions of them, are finite or countable
infinite (Theorem  \ref{srb}), derives the attention to try to
define and find sufficient conditions to reduce as much as
possible the set of observable measures.

Besides,  the \em reductions \em of the space of observable
measures will work as \em Generalized Ergodic Attractors, \em even
in the case that this reduction can not be done as much as to
obtain physical measures.
We first prove  that the reducibility of the set $\cal O$ of
observable measures for $f$ must be defined carefully, because in
the following sense, this set $\cal O$ is  minimal.
\em
\end{prueba}

\begin{theorem} \label{teoremaminimal} \em
{\bf (Reformulation of  Theorem \ref{toeremaminimal0})}

\em Let $f\colon M \mapsto M$  be any given continuous map in the
compact manifold $M$.

The set ${\cal O }_f $ of all its observable measures belongs to
the family
$$\aleph= \{ {\cal K} \subset {\cal P}: {\cal K} \mbox{ is compact
and } p\omega (x) \subset {\cal K} \mbox{ for Lebesgue almost
every point } x \in M \}$$

 Moreover $${\cal O}_f = \bigcap _{{\cal K} \in \aleph} {\cal
K}$$ and thus ${\cal O}_f$ is the unique minimal set in $\aleph$.
\end{theorem}

{\em Proof:} For simplicity let us denote ${\cal O}={\cal O}_f$.
Given any  subset $\cal K \subset \cal P$ (this $\cal K$ is neither
necessarily in $\aleph$ nor necessarily compact), let us
consider: \begin{equation} \label{2}
 A({\cal K})= \{x \in M: p\omega(x)\cap
{\cal K} \neq \emptyset \}, \;\;\;\;\;\;\; C({\cal K})= \{x \in M:
p\omega(x)\subset {\cal K} \}\end{equation}

 It is enough to prove  that $m(C({\cal O})) =1$ and that ${\cal K} \supset
\cal O$ for all ${\cal K} \in \aleph$. Let us first prove the
second assertion.

To prove that ${\cal K} \supset \cal O$ it is enough to show that
$\mu \not \in {\cal O}$ if $\mu \not \in {\cal K} \in \aleph$.

If $\mu \not \in {\cal K}$ take $\varepsilon = \dist (\mu, {\cal
K})>0$.  For all $x \in C({\cal K})$ the set $p\omega (x)\subset
{\cal K} $ is disjoint with the ball $B_{\varepsilon }(\mu)$. But
almost all Lebesgue point $x \in C({\cal K})$, because ${\cal K} \in
\aleph$. Therefore $p\omega (x) \cap B_{\varepsilon }(\mu) =
\emptyset$ Lebesgue a.e. This last assertion and  Definition
$\ref{definicionobservable}$ and paragraph \ref{definicion2} imply
that $\mu \not \in {\cal O}$, as wanted.

Now let us prove that $m(C({\cal O})) =1$, which is the key matter
of this theorem. We know $\cal O$ is compact and not empty. So, for
any $\mu \not \in {\cal O}$ it is defined the distance $\dist (\mu,
{\cal O})>0$.  Observe that the complement  ${\cal O}^c$ of $\cal O$
in $\cal P$ can be written as the increasing union of compacts sets
${\cal K}_n$ (not in the family $\aleph $) as
follows:\begin{equation} \label{3} {\cal O}^c = \bigcup
_{n=1}^{\infty} {\cal K}_n, \;\;\;\;\;\;\;\;{\cal K}_n = \{\mu \in
{\cal P}: \dist (\mu, {\cal O}) \geq 1/n \} \; \supset  \; {\cal
K}_{n+1}\end{equation} Let us take the sequence $A_n= A({\cal K}_n)$
of sets in $M$ defined in (\ref{2}) at the beginning of this proof,
and denote $A_{\infty}= A({\cal O}^c)$. We deduce from (\ref{2}) and
(\ref{3}) that:
$$A_{\infty } = \bigcup _{n=1}^{\infty} A_n, \;\;\;\; m(A_n)\rightarrow m(A_{\infty}), \;\;\;\;
A_{\infty} = A({\cal O}^c)$$ To finish the proof is thus enough to
show that $m(A_n)=0$ for all $n \in \mathbb{N}$.

In fact,  $A_n= A({\cal K}_n)$ and ${\cal K}_n$ is compact and
contained in ${\cal O}^c$. By Definition \ref{definicionobservable}
and paragraph \ref{definicion2} there exists a finite covering of
${\cal K}_n$ with open balls ${\cal B}_1, {\cal B}_2, \ldots, {\cal
B}_k$ such that \begin{equation}  \label{4}
 m(A({\cal B}_i))= 0 \;\;\;
\mbox{for all } i = 1, 2, \ldots , k\end{equation}
 By (\ref{2})  the
finite collection of sets $A({\cal B}_i); \; i= 1,2, \ldots, k$
cover $A_n$  and therefore (\ref{4}) implies $m(A_n) = 0$ ending the
proof. $\Box$

\vspace{.2cm}

As shown in the examples  of Section \ref{ejemplos2}, there exist
maps whose spaces $\cal O$ of observable measures are irreducible
and maps for which they are reducible. Also there exist   maps
that do not have
 irreducible subsets in  $\cal O$.
 In section \ref{chains} we define chains and
co-chains of reductions. Those systems having physical measures can be
characterized also according to the existence of adequate decreasing sequences (chains)
of generalized attractors.

 For further
uses we define:

\begin{definition} \em \label{definicionAttracSize}   {\bf
(Diameter and Attracting Size.)}
Let $\cal O$ be the set of the observable measures for $f$. Let
${\cal O}_1 $ be a reduction or generalized attractor of
$\cal O$.

  \em The  diameter of ${\cal O}_1$ \em  is
 $\max\{dist (\mu, \nu): \mu, \nu \in
 {\cal O}_1\}$.
\em The attracting size of ${\cal O}_1 $ \em is $m(B({\cal
O}_1))$, where $B({\cal O}_1)$ is the basin of attraction of
${\cal O}_1$ (see definition \ref{definicionbasin}).
\end{definition}

By definition of generalized attractor its attracting size is positive.
%If the purpose is to find the minimal generalized attractors (supposing they exist), or at least
%a partition of the space of observable measures ${\cal P}$ into different generalized attractors, such that their diameters in ${\cal P}$ are as small as possible,    but all of them have positive attracting size, we should
%classify the phase space $M$ according to which basin of the generalized attractor each point $x \in M$
%belongs to. Then we should try to know when a given compact subset
%of the space of observable measures is a reduction, i.e. a generalized attractor. That is why we
%state, for  a sake of completeness,  some
%results that characterize the reductions.
 If the basin of attraction $B$ of some compact subset of ${\cal
O}$ has positive Lebesgue measure, then there exists a compact set
$K \subset B $ with positive Lebesgue measure. Thus we obtain the
following  characterization of all the reductions of $\cal O$, as
a consequence of Egoroff Theorem:

\begin{proposition} \em \label{teoremaCU} {\bf (Generalized attractors and
uniform
convergence.)} \em

 The subspace ${\cal
O}_1$ is a reduction of the space $\cal O$ of the observable
measures for $f$, \em (i.e. ${\cal O}_1$ is a generalized attractor
for $f$), \em if and only if there exists a positive Lebesgue
measure set $K \subset M$ such that
$$\mbox{\em dist} \left(\frac{1}{n}\sum_{j=0}^{n-1}\delta _{f^j(x)}, \;
{\cal O}_1 \right) \rightarrow 0 \;\; \mbox {uniformly in } x \in
K $$  \em Note: The set $K \subset M$  is not necessarily
$f$-invariant.
\end{proposition}

{\em Proof: } Let us call $B$ to the basin of attraction of ${\cal
O}_1$ (see definition \ref{definicionbasin}). We have $m(B)>0$ and
therefore,  the sequence in \ref{teoremaCU}  converges to 0
$m$-a.e. $ x \in B$.
 The direct result is now a straightforward consequence of
 Egoroff Theorem and its converse is obvious. $\Box$

 The following is other
characterization of the reductions  of the space  of observable
measures for $f \in C^0(M)$, in terms of the invariant subsets in
$M$ that have positive Lebesgue measure:

\begin{proposition} \em \label{proposicionrestriccion} {\bf (Restricting the map
to reduce the set of observable measures.)}

\em  The subspace ${\cal O}_1$ is a reduction  of the space ${\cal
O }_f$ of the observable measures for $f$ (i.e. ${\cal O}_1$ is a
generalized attractor for $f$) if and only if ${\cal
O}_1={\cal O }_{f_1}$,    where ${\cal O}_{f_1}$ is the set of all
observable measures of the map $f_1 = f|_C$, obtained when $f$ is
restricted to  some invariant set $C \subset M$ that has positive
Lebesgue measure. Besides $C$ can be chosen as the basin attraction ${B({\cal
O}_1)}$ of ${\cal O}_1$.
\end{proposition}

{\em Proof: } This Theorem is a corollary of Theorem
\ref{teoremaminimal}. In fact, to prove the converse statement apply
\ref{teoremaminimal} to $f|_C$ instead of $f$, taking $C={B({\cal
O}_1)} $ where ${\cal O}_1$ is the given reduction of ${\cal O}_f$.
To prove the direct result apply also  \ref{teoremaminimal} to
$f|_C$, but now taking $C$ as the given invariant subset in $M$ with
positive Lebesgue measure. $\Box$

\begin{prueba}  \em \label{pruebaGenErgAttract}
{\bf -  Proof of Theorem \ref{toeremaminimal0}}:
 By Proposition
\ref{proposicionrestriccion} and Theorem \ref{teoremaminimal}
applied to $f |_ C$ where $C = B({\cal O}_1)$, we have that Lebesgue
almost all $x \in C$ verifies $p\omega(x) \subset {\cal O}_1$, and
any $\mu \in {\cal O}_1$ is observable for $f|_C$. Take any $\mu \in
{\cal O}_1 \setminus {\cal K}$. There exists an open ball
$B_{\varepsilon }(\mu )$ that does not intersect ${\cal K}$. As
$\mu$ is observable for $f|_C$, (see Definition
\ref{definicionobservable}), $p\omega(x)$ is not contained in ${\cal
K}$ for a set of $x \in C$ with positive Lebesgue measure. Therefore
$m (C \setminus B({\cal K})) >0$ as wanted. $\Box$
\end{prueba}

\section {Examples.}  \label{ejemplos2}

The  examples in this section are well known or very simple, but
give a  scenario of  the possible dynamics, in terms of the
statistics given by the time mean sequence (\ref{1}). In fact, they
are paradigmatic of some  different classes of continuous dynamical
systems  $f \in C^0(M)$.  After Theorem
\ref{teoremaDescomposicionGenErgAttr0}), some of these examples may
appear joint with the others, in each of the basins of the (up to
countable infinitely many) generalized attractors, in which the
complete set ${\cal O}_f$ of observable measures decomposes. In a
general case, the complete topological dynamics may be much more
complicated,
 since the basins of the infinitely many
generalized attractors of $f$, may be topologically   riddled   in $M$, as explained
in the discussion at the end of this section.
\begin{example} \em \label{ejemploPozo}
For a map with a single periodic point $x_0$, being a topological
sink whose topological basin is $M$ almost all point, the set
${\cal O}$ has a unique measure that is the $\delta$-Dirac measure
supported on $x_0$.
\end{example}
\begin{example} \em \label{ejemploAnosov}
For  any transitive Anosov $C^{1 + \alpha}$ diffeomorphism  the set
$\cal O$ is irreducible containing uniquely the SRB measure $\mu$.
But there are also infinitely many other ergodic and non ergodic
invariant probabilities, that are not observable (for instance the
equally distributed Dirac delta measures combination supported on a
periodic orbit). In particular, if
    $f $ preserves a probability $\mu$ equivalent to the Lebesgue measure, then $f$ has
$\mu$ as the unique observable, and thus physical, probability.  This result is generalized
also for some subclass of $C^{3}$ diffeomorphisms on the two-torus, conjugated to transitive Anosov,
   with a non
hyperbolic fixed point in which the derivative of $f$ has double
eigenvalue 1 and a single eigendirection (\cite{nosotros}). But is
is mostly open for other conjugated to transitive Anosov, even if
they are analytic.
\end{example}

\begin{example} \em \label{ejemploHu}
In \cite{huyoung} it is studied the class of $C^2$ diffeomorphisms
$f$ in the two-torus obtained from a transitive Anosov, and in the
same conjugation class of the Anosov, when the unstable eigenvalue
of a fixed point $x_0$ is weakened to be 1, maintaining its stable
eigenvalue strictly smaller than 1 and maintaining also the uniform
hyperbolicity in each iterate outside a neighborhood (non invariant)
of the fixed point. It is proved that $f$ has a single physical
measure that is the Dirac delta supported on $x_0$ and that its
basin has total Lebesgue measure. Therefore this is the single
observable measure for $f$, although there are infinitely many other
ergodic invariant measures. As the physical measure is supported in
a fixed point $x_0$, statistically $x_0$ acts as a sink, attracting
the sequences of time averages of Lebesgue almost all orbit.
Nevertheless, as $f$ is conjugated to Anosov, it is topologically
chaotic (i.e. expansive, or sensible to initial conditions).
\end{example}

\begin{example} \em \label{ejemploCao}
The diffeomorphism $f\colon [0,1]^2 \mapsto [0,1]^2; \;\; f(x,y) =
(x/2, y)$ has the set ${\cal O}$ of observable measures as the set
of  Dirac delta measures $\delta _{(0,y)}$ for all  $ y \in [0,1]$.
In this case $\cal O$ coincides with the set of all ergodic
invariant measures for $f$, it is infinitely reducible (i.e. $\cal
O$ is reducible and any reduction of ${\cal O}$ is also reducible).
Not all $f$-invariant measure $\mu$ for $f$ is observable: for
instance, the one-dimension Lebesgue measure on the interval
$[0]\times [0,1]$ is invariant and is not observable. This example
shows that the set ${\cal O}$ is not necessarily convex.

%-------------------------------------
%
%In this example, any measure being the convex combination of Dirac's delta is stochastically observable.
%
%-------------------------------------

\end{example}

\begin{example} \label{ejemploInfinitosPozos}
 \em
 The maps exhibiting infinitely many simultaneous hyperbolic sinks, constructed from
 Newhouse's theorem (\cite{newhouse}) has a space ${\cal O}$ of observable
 measures that is reducible. But it has infinitely many reductions
 (the Dirac delta supported on the hyperbolic sinks), each of them being
 irreducible. Also the maps exhibiting infinitely H\'{e}non-like
 attractors, constructed by Colli in \cite{colli},  has a space of
 observable measures that is reducible, having infinitely many
 reductions (the physical measures supported on the H\'{e}non-like
 attractors), each one that is irreducible.

\end{example}
\begin{example}
\em \label{ejemploFeigenbaum}
Consider  the quadratic family $\{f_t\}_t$ in the interval $I$, and in this family, the map $f  $
where the first cascade or period doubling bifurcating maps converge. It has a single attractor $A$
which is   Feigenbaum-Coullet-Tresser. This attractor $A$ is formed by a single orbit, whose closure
$K$ is a Cantor set having as extremes of its gaps, the  future orbit  of the critical point.
For all $x \in A$, $f^n(x)$   moves quasi-periodically in a single orbit
(with quasi-periods $2^n$ for all $n \geq 1$) and attracts topologically all the points of $I$,
except those of   countable many   periodic hyperbolic repellors (with periods   $2^n$, for all
$n \geq 0$). The map $f$ is infinitely doubling renormalizable and has a single observable
probability $\mu$ (and thus physical measure) supported on $K$. This physical measure
is constructed as follows:  For any fixed $n \geq 1$
call $\{K_{i,n}: \, {   0 \leq i \leq 2^n-1} \}$ to  the family of
$2^n$ atoms of generation $n$, i.e.   the sets $ K_{i,n} : \; 0 \leq i (mod \ 2^n) \leq 2^n-1 $
are the pairwise disjoint compact invervals
such that $f(K_{i,n}) \subset K_{i+1,n}$, and $K = \bigcap _{n \geq 1}
\bigcup_{i= 0}^{2^n-1} K_{i,n}$. Then define the probability $\mu$
such that $\mu(K_{i,n}) = 1/2^n$ for all
$n \geq 1$ and for all
$0 \leq i   \leq 2^n-1$.

\end{example}

\begin{example}
\em \label{grad} (Example 1 in \cite{araujo2}). Define
the function $\phi \colon [- \, \frac1\pi,\frac1\pi] \mapsto
\mathbb{R}$, $\phi(s)=s^4 \sin(\frac1s)$ and identify the extremes
of the interval to obtain $\varphi \colon S^1 \mapsto \mathbb{R}$.
Consider the time one of the gradient flow given by
$\dot{x}= \nabla \varphi(x)$. This map has infinitely many sources
and sinks, which accumulate at 0. The physical measures are the
infinitely many sinks. The Dirac  delta on the accumulation point of
the sinks and
sources, is not a physical measure, but it is an observable measure.
It is besides the unique stochastically stationary measure. We conclude
 that physical measures, even if their basins attract
Lebesgue a.e.,
 do not necessarily include the stochastically stationary probabilities but, at least
 in this example, observable measures  include them.
\end{example}

%\begin{example}
%\em \label{grad2} This example shows that observable measures and
%observable stationary measures can be different. Take the former
%example (\ref{grad}) and \lq\lq blow up" a source, in such a way
%that the diffeomorphism a at such interval be the identity. The set
%of observable stationary measures is the same one as before, but the
%set of observable measures is increased by the fixed points.
%\end{example}
%-----------------------------------------
%

\begin{example} \em \label{ejemplodossillas}
The following example, due to Bowen (see also  example 2 in
\cite{araujo2}), shows that the space of observable measures may be
formed by measures that are partial limits of the sequences of time
averages of the system states and that this sequence may be not
convergent for Lebesgue almost all points. Consider a diffeomorphism
$f$ in a ball of $\mathbb{R}^2$ with two hyperbolic saddle points
$A$ and $B$  such that the unstable global manifold $W^u(A)$ of $A$
is a embedded arc that coincides (except for $A$ and $B$) with the
stable global manifold $W^s(B)$ of $B$, and the unstable global
manifold $W^u(B)$ of $B$ is also an embedded arc that coincides with
the stable manifold $W^s(A)$ of $A$. Let us take $f$ such that there
exists a source $C$ in the open  ball $U$ with boundary $W^u(A)\cup
W^u(B) $, and all orbits in that ball $U$ have $\alpha$-limit  set
$C$ and $\omega $-limit set $W^u(A) \cup W^u(B)$. If the eigenvalues
of the derivative of $f$ at $A$ and $B$ are well chosen, then one
can get that the time average sequences of  the orbits in $U
\setminus \{C\}$ are not convergent,  have at least one subsequence
convergent to the Dirac delta $\delta_A$ on $A$ and  have other
subsequence convergent to the Dirac delta $\delta _B$ on $B$.

Due to Theorem \ref{teoremaconvexo}, for each $x \in U
\setminus\{C\}$ there are non countably many probability measures
which are the limit measures of the time average sequence of the
future orbit starting on $x$. All these measures are invariant
under $f$ and therefore, due to Poincar\'{e} Recurrence Theorem (see
\cite{mane}), all of them are supported on $\{A\} \cup \{B\}$. Due
to this last observation and due to Theorem \ref{teoremaconvexo}
all the convex combinations of $\delta _A$ and $\delta _B$ are
limit measures of the sequence of time averages of any orbit
starting at $U \setminus \{C\}$ and conversely.

Therefore the set ${\cal O}$ of observable measures for $f$
coincides with the set of convex combinations of $\delta _A$ and
$\delta _B$. The set ${\cal O}$ is irreducible and formed by
non-countable many probability measures. It is not  the set of all
invariant measures; in fact the measure $\delta _C$ is not
observable.

This example also shows that the observable measures are not
necessarily ergodic.

A different exact adjustment in the eigenvalues of the two saddles
$A$ and $B$ allows a different example, for which all the convergent
subsequences of (\ref{1}) converge to the same exact previously
chosen convex combination $\mu$ of $\delta_A$ and $\delta_B$. So,
there is a physical measure that attracts   the time averages of the
orbits of $U$ and that is not ergodic. This proves that physical
measures are not necessarily ergodic.

In \cite{araujo2} it is shown that this physical measure, which is a
  convex combination $\mu$ of $\delta_A$ and
$\delta_B$, is stochastically stable, even being non ergodic.

\end{example}

\begin{example} \label{ejemploExpanding} \em
The $C^{1 + \varepsilon}$ expanding maps $f\colon S^1 \mapsto S^1$
in the circle (i.e. $f'(x) >1 \; \forall \; x \in S^1$), have been
extensively studied, they present a unique SRB and physical measure
attracting Lebesgue a.e. In \cite{cq} it is shown that also in the
$C^1$ topology, generically $f$ has a unique physical measure.
Nevertheless this physical measure is not SRB (it is not absolute
continuous respect to Lebesgue) in those generic $C^1$ (not $C^{1+
\varepsilon})$ examples. They  show   that, for   $C^1$ uniformly
hyperbolic maps, (that are besides topologically mixing), if there
is a unique observable   (and thus physical) probability, this
measure is not necessarily SRB.

On the other hand,   there exists $C^1$  examples of expanding maps
for which a non ergodic invariant measure $\mu$ is equivalent to
Lebesgue (\cite{Q3}). Thus, %the observable measures  in those
%examples, are  ergodic components of $\mu$ (see Remark
%\ref{teoremaLebesgue0}) $\clubsuit$ ?`por qu\'{e} las componentes
%erg\'{o}dicas eran observables? (olvid\'{e} el razonamiento):
there is not a unique observable measure. This shows that, out of
the $C^{1 + \varepsilon}$ case, when uniform hyperbolicity and
topological mixing hold,   the set of observable measures is not
necessarily reduced to a single probability, even if the sequence
(\ref{1}) converges Lebesgue a.e.

The case of $f \in C^0$ stands in contrast with the former
situation. In \cite{misiurewicz} (see Theorem~3.4) is proved in
particular that there exist maps $f\colon S^1 \to S^1$ topologically
conjugate with $g_d\colon S^1 \to S^1$, $g_d(x)=dx$ such that for
Lebesgue almost every point $x$ in $S^1$ and every $f$-invariant
measure $\mu $, some subsequence of the sequence (\ref{1}) converges
to $\mu$.

% In Theorem \ref{teoremon} we prove that the  observable measures
% of the $C^1$ expanding maps $f$ are contained in the set
% of the equilibrium states of $-\log f'$.

 %Therefore, in particular a $f$-invariant measure
% $\mu$ is observable  only if the Ruelle inequality is an equality (PLY % equality), i.e. the   entropy  of $\mu$ is equal to $\int  \log f' \; d \mu$.
% as
 \end{example}

  \vspace{.3cm}
  {\bf Discussion: } When formulating the theorems about chains and
  co-chains of generalized attractors in Section \ref{seccionphysical},
we have in mind    the three  paradigmatic    different statistical
 dynamical behaviors that $C^0$ systems
 may exhibit (in relation to the limit probabilities of the sequence (\ref{1})):

 First,
 in Examples \ref{ejemploCao} %and $\clubsuit$ corregir esto, \ref{ejemploExpandingNoErgodico}
   Lebesgue a.e. orbit defines a convergent subsequence  (\ref{1}) of probabilities,
   but no  physical measure exist (at least
 for a subset of initial states with positive Lebesgue measure).

 Second,
 in Examples \ref{ejemplodossillas} %and  $\clubsuit$ corregir esto, \ref{ejemploExpandingNoErgodico}
  Lebesgue a.e. orbit defines a non convergent subsequence (\ref{1}).

Third, in Examples \ref{ejemploExpanding} (the $C^2$ or the $C^1$
generic expanding maps), \ref{ejemploFeigenbaum} (the Feigenbaum
attractor),   \ref{ejemploInfinitosPozos} (infinitely many
coexisting sinks or H\'{e}non like attractors), \ref{ejemploAnosov} and
\ref{ejemploHu} ($C^2$ conjugated to Anosov in the two-torus),
Lebesgue almost all orbit defines a convergent subsequence (\ref{1})
that converges to a physical measure, and there are at most
countable many such physical measures.

 To state the Theorems of Section \ref{seccionphysical},  we are thinking in those   different
 statistical dynamical behaviors, but
 not as dynamical systems that are topological isolated one from the others.
Precisely, we are considering that the  basins of the  generalized
attractors of all those examples may be  topologically immersed in
a larger dimension compact manifold $\widehat M$, in
such a way that  become mutually riddled  (i.e. they are dense
subsets of $\widehat M$).

\section{Equilibrium states and observable measures for $C^1$ expanding maps.}
We will develop the  theorems in
this section  for order preserving expanding maps in the circle $S^1$, but the proofs also work
   for order reversing
expanding maps in   $S^1$. Some of the results also
work for expanding maps, similarly defined, in appropiate manifolds of dimension larger than one.

\begin{definition} \em
A $C^1$ (order preserving) map $f: S^1 \mapsto S^1$ is expanding if $f'(x) >1$ for all $x \in S^1$.
(If $f$ is order reversing then it is expanding if $-f'(x) >1$.)

In particular, for all integer $d >1$ we denote $g_d : S^1 \mapsto S^1$ to the linear expanding map
$$g_d(x) = dx \; \forall \, x \in S^1$$ If $f$ is expanding, then
there exist a unique integer $d >1$, called the degree of $f$, such
that $f$ is topologically conjugated to $g_d$.  Thus, all the
topological properties of $g_d$ are translated to $f$.

We denote ${\mathcal E}^1 \subset C^1 (S^1)$ to the family of $C^1$
expanding maps in the circle $S^1$. For all $r \geq 1$, we denote
${\mathcal E}^r \subset {\mathcal{E}^1}$ to the $C^r$ expanding maps
in $S^1$. Denote $\mbox{Homeo}(S^1)$ to the set of homeomorphisms on
the circle $S^1$ with the $C^0$ topology. Finally, for all $f \in
{\mathcal{E}^1}$,  denote $\mbox{Conj}(f) \subset \mbox{Homeo}(S^1)
$ to the set of all the conjugacies between $f$ and $g_d$.

\end{definition}
For a seek of completeness we state here a Lemma  from \cite{cq}.
For each $x \in S^1$ denote $U_x$ to the $C^1 $ open and dense set
 of all expanding maps $f$ in $S^1$ such that $f(x) \neq x.$
\begin{lemma} \em
{\bf   }
 \label{conj}
For each $x \in S^1$ there is a continuous map
$\Pi_x\colon U_x \to \mbox{Homeo}(S^1)$ such that $\Pi_x(f) \in \mbox{Conj} (f)$ for each $f$.
In particular, given $f \in \mathcal{E}^1$ of degree $d$, there is a neighborhood $U$
of $f$ on which there is a continuous choice of conjugacies to the map $g_d$.

\end{lemma}
{\em Proof: } See Lemma 1 of  \cite{cq} and Theorem 2.4.6. of
\cite{kh}.

\begin{definition} {\bf Pressure and equilibrium states.} \em  \label{definicionPresure}

Let
  $f   $ a $C^1$- expanding map in the circle. We denote $\psi = -\log f' \in C^0 (S^1, \mathbb{R})$
and ${\cal P}_f$ to the set of $f-$invariant Borel probabilities in $S^1$. For $\nu \in {\cal P}_f$
we denote $h_{\nu}(f)$ to the entropy of $\nu$.

Let $\varphi \in C^0(S^1, \mathbb{R})$. The pressure of $\varphi $,   is
$$P_f(\varphi)=\sup_{\nu\in \mathcal{P}_f} \left\{h_\nu(f) + \int \varphi \, d\nu\right\}.$$

A measure $\mu \in \mathcal{P}_f$ is \em an equilibrium state for $\varphi$  \em if $$h_\mu(f)
+ \int \varphi \, d\mu = P_f(\varphi).$$
In particular $\mu \in {\cal P}_f$ is \em an equilibrium state for $\psi = - \log f'$ \em if
$$h_\mu(f)
=  \int \log f'(x) \, d\mu(x) + P_f(-\log f').$$
We denote as $ES_f \subset {\cal P}_f$ to the set of $f$-invariant probabilities that are
equilibrium state for $\psi = - \log f'$.
\end{definition}
For a seek of completeness we recall  well known results in the following Theorem \ref{teoremaRuelle}
and Corollary \ref{corolarioRuelle}:
\
\begin{theorem} \label{teoremaRuelle} {\bf   Ruelle inequality and Pesin-Ledrappier-Young Equality
for $C^{1 + \alpha}$ expanding maps.}
If $f \in {\mathcal E}^{1 + \alpha}$ \em (i.e. $f$ is a $C^1$ expanding map in the circle such that
$f'$ is $\alpha >0$-H\"{o}lder continuous),  \em  then  there exists a unique $f$-invariant and ergodic
measure $\mu$ that is the equilibrium state for $\psi = - \log f'$. Besides $\mu$ is the unique
$f$-invariant measure that is absolute continuous respect to  the Lebesgue measure $m$,
 and the pressure   $P_f(\log f') = 0$. In other words, for all $\nu \in {\cal P}_f$ the inequality
 of Ruelle holds:
 $$h _{\nu}(f) \leq \int \log f'(x) \; d \nu.     $$
 Besides, this last is an equality if and only if $\nu \ll m$, and this holds if and only if $\nu = \mu$.
\end{theorem}
 {\em Proof: }  See Theorem  6.3.8 of \cite {ke}.

 \begin{corollary}
 \label{corolarioRuelle} {\bf Ruelle inequality for $C^1$ expanding maps.}
If $f \in {\mathcal E}^{1 }$ \em (i.e. $f$ is a $C^1$ expanding map in the circle),
 \em  then   the pressure   $P_f(\log f') = 0$.
 In other words, for all $\nu \in {\cal P}_f$ the inequality
 of Ruelle holds:
 $$h _{\nu}(f) \leq \int \log f'(x) \; d \nu .     $$
 \end{corollary}

 {\em Proof: } We reproduce the proof of Lemma 2 in \cite{cq}.
   After the PLY Equality of Theorem \ref{teoremaRuelle} it
   follows that $$P_f (-\log f') = h_{\mu} (f) - \int \log f' \; d \mu =  0 \ \mbox{ if } \
   f \in \mathcal{E}^2$$ for the unique
   $f$-invariant measure $\mu $ which is absolute continuous
    respect to the Lebesgue measure $m$.

 If   $f \in \mathcal{E}^1$ has degree   $d$, after the lemma \ref{conj}, there exists  a
 neighborhood $U \subset \mathcal{E}^1$ of $f$, and for all $g \in U$,
  conjugacies $\gamma_g \in \mbox{Conj}(g)$, between $g \in U$
 and the linear expanding map $g_d$ of degree $d$,
  such that the application
 $g \to \gamma_g$ is continuous on $U$. Denote $\psi_{g} (x) = - \log g' (x)$ for all $g \in U$.
 Take $\{f_i\}_{i\in \mathbb{N} }
 \subset \mathcal{E}^2$ such that $f_i \to f \in \mathcal{E}^1$ (the convergence is in the $C^1$ topology).
  Then $\psi_{f_i}\circ \gamma_{f_i}
 \to \psi_f\circ \gamma_f \in C^0(S^1, \mathbb{R})$. Note that the pressure  $p_g(\varphi)$, for any
 fixed $g \in {\cal E}^1$,
 depends continuously on $\varphi \in C^0(S^1, \mathbb{R})$. Then $$P_{g_d}
 (\psi_{f_i} \circ \gamma_{f_i}) = P_{g_d} (\psi _f \circ \gamma_f).$$

Also note that, if $\gamma \in \mbox{Conj}(g)$ for some $g \in
{\mathbb{E}^1}$ and if $\varphi \in C^0(S_1, \mathbb{R})$, then
$$P_{g_d} (\varphi \circ \gamma)
 = P_g(\varphi).$$ Indeed, $\gamma$  induces a bijection $\gamma^*$ between the $g_d$
 invariant measures
 $\nu \in {\cal P}_{g_d}$, and the  $g$ invariant measures $\gamma^* \nu \in {\cal P}_g$.
  Then $\int \varphi \circ \gamma \; d \nu =
 \int \varphi \; d \gamma^* \nu$ for all $\nu \in {\cal P}_{g_d}$. Besides, since $\gamma^*$
 is a measure-theoretic isomorphism, the entropies coincide $h_{\nu}(g_d)= h_{\gamma^* \nu} (g)$.

We conclude that  $$0=P_{f_i} (\psi_{f_i}) = P_{g_d}(\psi_{f_i}\circ
\gamma_{f_i})=P_{g_d} (\psi_{f }\circ \gamma_{f })
   = P_f (\psi_{f }). \  \; \; \Box$$

 \begin{remark} \em

 If $f $ is a $C^1$ expanding map in the circle that is not $C^{1 + \alpha}$, then
 (except the inequality of Ruelle), the  thesis of Theorem \ref{teoremaRuelle}
 does not necessarily hold. In fact, the uniqueness of the equilibrium state of $-\log f'$,
 or its absolute continuity
 respect to Lebesgue, may fail, as  in the following cases:

 On one hand,
 the PLY equality may hold for a unique equilibrium state that is a
physical measure  but singular   respect to Lebesgue (\cite{cq}).

On the other hand, the PLY equality may hold for some invariant probability that is
 absolute continuous
respect to Lebesgue, but not ergodic: see \cite{Q3} and Theorem of Ledrappier (Theorem 2 in \cite{walters}). Therefore  its ergodic components are also equilibrium
states of $-\log f'$.

 \end{remark}
 The following theorem states that any observable measure is an
 equilibrium state. It  is a stronger version of Theorem 6.1.8 of the book of Keller \cite{ke}{\footnote {Our definition
of observability is weaker than the definition in \cite{ke}}}, but its proof is in essence the same.

\begin{theorem} \label{teoremon}
{\bf }

Let $f$ be a $C^1$ expanding map on the circle $S^1$. Then, any
observable measure of  $f$ and any convex combination of observable
measures of $f$, is an equilibrium state for $\psi= - \log f'$,  and
its pressure is equal to 0, i.e. any observable measure $\mu$ of $f
$ satisfies the PLY equality for the entropy:
$$h_\mu(f)= \int \log f' \, d\mu.$$

\end{theorem}

We prove Theorem \ref{teoremon} in the subsection
\ref{pruebadirectoTeoremon}. Let us now state its
 Corollaries. The first Corollary is a well known result. Nevertheless a new point of view for its
  proof
   rises
 from Theorems \ref{teoremon},  \ref{teoremaexistencia0} and \ref{toeremaminimal0}.

 \begin{corollary} \label{corolarioAgregado}
 A $C^1$ expanding map on the circle always has a non empty    set $ES_f$
 of probability
 measures that verify the PLY equality of the entropy. Besides,
 $f$ has a physical measure $\mu$ attracting Lebesgue a.e. if and only
 if the observable measure is unique, and this happens if there is a unique probability $\mu \in ES_f$.
 \end{corollary}

 {\em Proof: } After Theorem  \ref{teoremaexistencia0}: ${\cal O}_f \neq \emptyset$.
 After Theorem \ref{teoremon}, the closed convex hull of ${\cal O}_f$ is contained in $ ES_f$, so $ES_f \neq \emptyset$.
 Finally, from Theorem \ref{toeremaminimal0}, $f$ has a physical measure $\mu$ attracting
 Lebesgue a.e. if and only if ${\cal O}_f = \{\mu\}$ and this trivially holds
 if $ES_f = \{\mu\}  \ \ \Box$

  We say that a probability measure is atomic if it is supported on a finite set.

\begin{corollary}  \label{corolarioteoremon}
There is no atomic ergodic observable measure of a $C^1$ expanding map.
\end{corollary}
{\em Proof of the corollary \ref{corolarioteoremon}:}  By
contradiction, if $\mu $ is an atomic ergodic observable measure of
an expanding map $f$, then $h_\mu(f)=0$. As $ \psi = -\log f'   <0$,
and $\mu $ is an equilibrium state for $\psi$, then the pressure
$P_f(\psi)$ is negative, contradicting Corollary
\ref{corolarioRuelle}. $\ \Box$

In the following definitions and lemmas, we reproduce those of the book of Keller (\cite{ke}) about
the equilibrium states  in its section \textsection 4.4,
 applied in particular to expansive $C^1$ maps $f$ in the circle $S^1$.
Recall that if $f \in {\cal E}^1$, the set of equilibrium states of $\psi= - \log f'$ is $ES_f \subset {\cal P}_f$, i.e.
$\mu \in ES_f$ if and only if the PLY equality
of the entropy holds: $h_{\mu}(f) = \int  \log f' \; d \mu$.

\begin{prueba}
{\bf Notation:}  \em \label{notacionKsubr}
 For
all $\nu \in {\cal P}_f$  we denote $$V_f (\nu) = h_{\nu}(f) -\int
\log f'(x) \; d \nu(x) $$ Due to Ruelle inequality (Corollary
\ref{corolarioRuelle}):
$$V_f (\nu ) \leq 0 \; \; \; \forall \nu \in {\cal P}_f
. $$
For all $r \geq 0$ we denote $${\cal K}_r = \{\nu \in {\cal P}_f: \;
V_f (\nu) \geq - r\}$$
In particular ${\cal K}_0 = ES_f$. For all $r \geq 0$
the set ${\cal K}_r$ is non empty, compact (in the weak$^*$ topology) and convex (join
the proof of
Theorem
4.2.3  of the book in \cite{ke},    with  Theorem 4.2.4 and Remark  6.1.10
of the same book).
\end{prueba}

For all integer $n \geq 1$  and all $x \in S^1$ denote
$\sigma_{n,x}$ to the (non necessarily $f$ invariant) probability of
the sequence (\ref{1}), called the \em empirical distribution of the
future orbit of $x$ up to time $n$, \em i.e.:
$$\sigma_{n,x} = \frac{1}{n} \sum_{j= 0}^{n-1} \delta_{f^j(x)}$$
where $\delta_x$ is the Dirac-delta probability supported on the point $x$.

As stated from the beginning of this paper, for $x \in S_1$ fixed,
we denote $$pw(x) \subset {\cal P}_f$$ to the set
of all probabilities that are the weak$^*$-partial limits (limits of the convergent subsequences)
 of the sequence  $\{\sigma_{n,x}\}_{n \geq 1}$ of the empirical distributions of the future orbit
 of $x$.

\begin{lemma}
\label{lemma1teoremon}  Let $f$ be a $C^1$ expanding map of the
circle $S^1$. Let ${\cal K}_r  $ the compact set of $f$ invariant
probabilities defined in \ref{notacionKsubr}. For all $r \geq 0$ and
for all open neighborhood ${\cal V} $ of ${\cal K}_r$ in the space
${\cal P}$ of all the (not necessarily $f$-invariant) probabilities,
the following inequality holds:
 $$\limsup_{n \rightarrow + \infty} \frac{1}{n} \log m \{x \in S_1: \sigma_{n,x} \not \in {\cal V}\}
  \; \leq  \; - r,$$
  where $m$ is the Lebesgue probability in the circle $S^1$.
\end{lemma}
{\em Proof:} This Lemma is the Proposition 6.1.11 of the book of
Keller \cite{ke}. All the hypothesis of that proposition\footnote{
The needed   hypothesis are restated according to
\lq\lq Corrections to Equilibrium states in ergodic theory"\ in
http://www.mi.uni-erlangen.de/~keller/publications/equibook-corrections.pdf
 published  by the author of the book.} hold in the case
that $f$ is a $C^1$ expanding map of the circle $S^1$: see Remark
6.1.10 of \cite{ke}.

\begin{prueba}

\label{pruebadirectoTeoremon}
{\bf Proof of Theorem \ref{teoremon}:} \em
As $ES_f$ is convex, it is enough to prove that
 ${\cal O}_f\subset ES_f $.  Consider, for any $r >0$ the   compact set
$K_r \subset {\cal P}_f$ defined in \ref{notacionKsubr}.
By definition, the   intersection of the (decreasing with $r$) family $\{K_r\}_{r >0}$
is the non empty compact set $${\cal K}_0  = \bigcap_{  n=1}^\infty {\cal K}_{\frac1n}$$

It is enough to prove that its basin of attraction $${B({\cal K}_0)}
= \{x \in S_1: pw(x) \subset {\cal K}_0\}$$ has total Lebesgue
measure. In fact, if we prove that $m(B({\cal K}_0)) = 1$, then
applying Theorem \ref{toeremaminimal0}, ${\cal O}_f \subset {\cal
K}_0 = ES_f$ as wanted.

 Now, we reproduce the proof of the part a. of Theorem
6.1.8 of \cite{ke}. Let $r >0$. We fix $0<\varepsilon < r $.  After Lemma \ref{lemma1teoremon}, for any
 weak$^*$ neighborhood ${\cal N}$
 of ${\cal K}_r \subset{\cal P}$, there exists $n_0$
 such that for any $n>n_0$: $m \{x:\sigma_{n,x} \in \mathcal{P} \setminus {\cal N}\}\leq
 e^{-n( r-\varepsilon)}.$  This implies that $\sum_{n=1}^\infty m (x:\sigma_{n,x}
 \in \mathcal{P} \setminus {\cal N})< +\infty.$  After the Borel-Cantelli Lemma it follows that
 $$m\left(\bigcap_{n_0=1}^\infty \bigcup_{n=n_0}^\infty \{x:\sigma_{n,x} \in \mathcal{P}
 \setminus {\cal N}\}\right)=0.$$ In other words,  for
 each open neighborhood ${\cal N}$ of ${\cal K}_r$, for $m$-a.e. $x \in S^1$  there
 exists $n_0 \geq 1 $ such that $\varepsilon_{n,x}\in {\cal N}$ for all $n\geq n_0$. Hence,
 $pw(x) \subset {\cal K}
 _r$ for
 $m$-a.e. $x \in S^1$. It follows that $pw(x) \subset \bigcap_{1 \leq n\in \mathbb{N}} {\cal K}_{1/n}
 = {\cal K}_0 =
 \{ \mu \in \mathcal{P}_f: V_f(\mu) = 0\} $ for $m$-a.e.$x \in S^1$  as wanted.
$\Box$

\end{prueba}

\section{Cardinality and decomposition of ${\cal O}$.}
\label{seccionphysical}

The first aim of this section is to state some results that
characterize the maps having physical measures whose basins attract
Lebesgue a.e., in terms of the cardinality of the set of its
observable measures. In particular we
prove  the part (b) of Theorem \ref{teoremaexistencia0} and
Theorem \ref{teoremaDescomposicionGenErgAttr0}, that were delayed up
to this section.

The second aim of this section,  for a seek of
completeness, is to cover in the theory,
 all the possible cases in $C^0(M)$, including those
maps that are singular respect to Lebesgue. To do that we
analyze   the chains and co-chains of  reductions of the space of
observable measures, even if no physical measures exists, or if one
or more than one physical measure exists, but the union of their
basins of attraction do not cover Lebesgue~a.e.

\begin{theorem} \em \label{teoremaCardinalidad} {\bf Cardinality of ${\cal O}$ and
physical measures.}

 Let $f\colon M \mapsto M$ be any continuous map in the compact
manifold $M$.

\em If the set $\cal O$ of the observable measures for $f$, or
some proper reduction ${\cal O}_1$ of $\cal O$, is finite or
countable infinite then there exists  physical measures $\mu$ for $f$, precisely of
$f|_{B({\cal O}_1)}$, where $B({\cal O}_1)$ denotes the basin of attraction in $M$ of ${\cal O}_1$.

Conversely, if there exists a physical measure $\mu$ for $f$ then, either the
space $\cal O$ has a single probability measure $\mu$, or it is reducible and there
exists
 a proper reduction  ${\cal O}_2$  of $\cal O$ with a single element.
\end{theorem}

{\em Proof:} The converse assertion is immediate. In fact, if $\mu
$ is physical then the basin of attraction of $\mu$ has positive
Lebesgue measure, and thus $\{\mu\}$ is a trivial reduction of
${\cal O}$. It is either a proper reduction or not. If not then
${\cal O}=\{\mu\}$.

Let us prove the direct assertion. Denote extensively as $\{\mu
_n: n \in N\}$ the finite or countable infinite reduction ${\cal
O}_1$ (proper or not) that is given in the hypothesis. (If it has
only a finite cardinality, then repeat one or more of its elements
in the extensive notation, but include all of them at least once).

By Proposition \ref{proposicionrestriccion} the space ${\cal O}_1$
is the set of all observable measures for the restriction  $f|_C$
of $f$ to some forward invariant set   with positive Lebesgue measure,
say $m(C)
>0$. Thus $C \subset B({\cal O}_1)$ and it is not restrictive to assume that
$C = B({\cal O}_1)$.

Rename if necessary $f|_C$ as $f$, ${\cal O}_1$ as $\cal O$, and
rename as $m$  the Lebesgue measure in $C $, (i.e: the restriction
to $C$ of the Lebesgue measure in $M$, which is then renormalized
to be a probability measure in  $C$). Resuming:
 \begin{equation}
 \label{[7]}
\mbox{For every  } x \in C: p\omega (x ) \subset {\cal O}=\{\mu_n: n
\in \mathbb{N}
 \} \end{equation}
Let us define  $C_n \in C$ to be  candidates of the basins of
attraction for the measures $\mu _n$, and relate their respective
Lebesgue measures $m(C_n)$ as follows:
\begin{equation}
\label{[8]}C_n = \{x \in C: \mu_n \in p\omega(x)\}; \;\;\;\;\;  C=
\bigcup_{n=1}^{\infty} C_n; \;\;\;\;\; \sum_{n=1}^{\infty}m (C_n)
\geq m(C) = 1
\end{equation} So $m(C_n)
>0$ for some $n \in \mathbb{N}$.

 To end the proof we shall show that
for  all $x \in C_n: \{\mu_n\}= p\omega (x)$. Due to (\ref{[7]}) and
(\ref {[8]}), it is enough to prove that $C_n \cap C_k = \emptyset$ if $\mu_n
\neq \mu_k$.

By contradiction, suppose that for some $\mu_n \neq \mu_k$ there
exists a point $x \in C_n \cap C_k \subset C$. Then, from (\ref {[8]}) we
have $\mu_n , \mu_k \in p\omega (x)$. Now we apply Theorem
\ref{teoremaconvexo}  and (\ref{[7]}) to conclude that the space
$\cal O$ is non-countably infinite. $\Box$

\begin{prueba} \em \label{srb} \label{corolarioCo-cadenaphysical}
{\bf Proof of Theorem \ref{teoremaexistencia0} (b)} It is straightforward consequence
of Theorems \ref{toeremaminimal0} and \ref{teoremaCardinalidad}. $\Box$
\end{prueba}

Now, let us analyze in an abstract theory, the complementary case: the set of all observable measures is
a non countable infinite compact subset of the space of invariant probabilities.

\begin{definition} {\bf (Chains of  reductions.)} \label{chains}
A chain of reductions  \em of the space ${\cal O}$ of the
observable
 measures for $f$ is a (finite or countable infinite) sequence
  $\{{\cal O}_n \} _{ n \in I \subset \mathbb{N}}$ of
 reductions or \em generalized attractors \em
 (see Definition \ref{definicionreduccion0})
 such that ${\cal O}_{i} \varsubsetneq{\cal
 O}_j$ if $i > j$ in the set  $I$ of natural indexes.

 We call \em length of the chain \em to its finite or countable
 infinite cardinality $\# I$.

Recall that, by definition of reductions (i.e. generalized
attractors) each ${{\cal O}_n}$ is a compact part of the set of the
observable measures.

\end{definition}

 For any chain $\{{\cal O}_n  \}_{n \in I \subset  \mathbb{N}}$ of reductions of the space
 of  observable measures for $f$, let $$d_n = \mbox{diam}
  ({\cal O}_n); \;\;\;\;\;\;\;\;\;\; s_n = \mbox{attrSize}({\cal O}_n)$$
  where $\mbox{diam}$ and attrSize denote respectively the diameter and the attracting
  size, defined in \ref{definicionAttracSize}.

  Observe that $d_n$ and $s_n$ are non negative decreasing sequences.
We denote    $\underline{d} = \lim d_n$ and $\underline s = \lim s_n$.

\begin{theorem} \em {\bf (Physical measures and chains.)} \label{teoremaPhysicalMeasuresAndChains}
 \em There exists a physical measure for $f$ if and only if there
exists a chain ${\cal O }_n$ of reductions of the space $\cal O$
of the observable measures of $f$, such that the sequence of its
diameters converges to zero and the sequence of its attracting
sizes converges to some $\alpha
>0$.

\end{theorem}
{\em Proof: } The converse statement is immediate  defining the
length 1-chain $\{\{\mu\}\}$, where $\mu$ is the given physical
measure. The direct result is also immediate if the  length of the
given  chain is finite. Let us see now the case when the chain
${\cal O}_n; \; n \in \mathbb{N}$ is infinite. As the sequence of
its diameters converges to zero, then $\cap _{n \in \mathbb{N}}
{\cal O}_n = \{\mu \}$ for some $\mu$. It is enough to show that the
attracting size  of $\mu$ is positive. Note that from the
construction of such $\mu$ we have that $C(\{\mu\})= \cap _{n \in
\mathbb{N}} C_n $ where $C_n$ denotes the basins of attractions
$C({\cal O}_n)$. These basins are a countably infinite decreasing
family of sets in $M$ with positive Lebesgue measures $s_n$.
Therefore, $attrSize (\{\mu\}) = m (C(\{\mu\}))= \lim s_n = \alpha
>0$, as wanted. $\Box$

\begin{definition} {\bf (Independence of generalized attractors and chains.)}
\em We say that two generalized attractors or reductions
of the space of observable measures \em are   independent   if the
basin of attraction of their intersection has zero Lebesgue
measure. \em

We note from Definition \ref{definicionbasin} that the basin of
attractions of two reductions ${\cal O}_1$ and ${\cal O}_2$
intersect exactly in the basin of attraction of ${\cal O}_1 \cap
{\cal O}_2$. Therefore:

\em Two ergodic attractors are independent if and only if  the
intersection of their basins has zero Lebesgue measure.\em

We say that \em two chains of reductions are independent \em if
each one of the chains  has a reduction that is independent with
\em some \em reduction of the other chain.

\end{definition}

\begin{definition} {\bf (Co-chains of reductions.)} \label{teoremaCadena}
A co-chain of reductions  \em of the space ${\cal O}$ of the
observable
 measures for $f$ is a (finite or countable infinite) family
  $\{{\cal O}_n; \; n \in I \subset \mathbb{N} \}$ of
 reductions or \em generalized attractors \em
 (see Definition \ref{definicionreduccion0})  \em that are
 pairwise independent. \em

 We call \em length of the co-chain \em to its finite or countable
 infinite cardinality $\# I$.

 {\bf Remark:} If the space $\cal O$ of all the observable measures
 for $f$ is irreducible, then $\{{\cal O}\}$ is the unique chain of
 reductions and also the unique co-chain.

\end{definition}

Now we state a slightly generalized version of a  known result in
the theory of Discrete Mathematics, the Theorem of  Dilworth
(\cite{liu}), applied to the chains and co-chains of generalized attractors:

\begin{theorem}{\bf (Reformulation of  Dilworth Theorem.)} \label{teoremaDilworth}
For any continuous map $f$ the supreme $k$ of the lengths of the
co-chains of reductions in the space of observable measures for
$f$, is equal to the supreme $h$ of  the number of pairwise
independent chains.

Moreover: for any co-chain of length ${ l}$ there is a family of $\  { l}$
pairwise independent chains, and conversely.
\end{theorem}

{\em Proof: } Any co-chain $\{{\cal O}_j, \; j\in J\}$ with length
 $\; l = \# J$ can be seen as a collection $ \{\{{\cal O}_j\}, \;  j
\in J \} $ of $l$ pairwise independent chains $P_j = \{{\cal
O}_j\}$, each chain $P_j$ with length one. So $k \leq h$.

Conversely, given any collection $\{ P_j, \;  j \in J\}$ of
pairwise independent chains, take a reduction ${\cal O}_1 \in P_1$
independent to some $\widehat{\cal O}_2 \in P_2$, and take $\breve
{\cal O}_2 \in P_2$ independent with $\widehat{\cal O}_3 \in P_3$.
As both reductions $\widehat{\cal O}_2 $ and $ \breve{\cal O}_2$
belong to the chain $P_2$,  one of them must be contained in the
other; thus their intersection, say ${\cal O}_2$, is also a
reduction of the chain $P_2$. Besides ${\cal O}_2$ is independent
with ${\cal O}_1 \in P_1$ and with $\widehat {\cal O} _3 \in P_3$.
Analogously  construct by induction a (finite or infinite)
sequence  $\{{\cal O}_j: j \in J\}$ of pairwise independent
reductions, such that  ${\cal O}_j \in P_j$. This sequence of
reductions is by definition a co-chain. Therefore, $h \leq k$.
 $\Box$

When no physical measure exist, or when some  of them exist but
their basins do not cover Lebesgue almost every point of the phase
space $M$, we will still state a equivalent condition for the space
being partitioned in (up to countably many) irreducible generalized
attractors, whose basins cover Lebesgue all orbit.

\begin{theorem} \label{teoremaCo-cadena}
\em {\bf (Co-Chains and irreducible attractors.)}
\begin{enumerate}
\item \em A map $f\colon M \mapsto M$ has \em(up to countable
infinitely many) \em irreducible attractors whose basins of
attraction cover Lebesgue almost all point in $M$ if and only if
there exist a co-chain of reductions of the space $\cal O$ of the
observable measures for $f$ such that  the \em (finite or countable
infinite) \em sum of its attracting sizes is 1.

In this case:

\item   The irreducible attractors are all physical measures if
and only if the diameter of the reductions are all zero.

\item  The \em (finite or countable infinite) \em number of such
irreducible generalized attractors is equal to the supreme $k$ of
the lengths  of the co-chains of reductions for $f$ and to the
supreme $h$ of the number of independent chains for $f$.

\end{enumerate}
\end{theorem}

{\em Proof: }  From Definition
$\ref{definicionAttracSize}$ we obtain that $f$ has generalized attractors
whose basins cover Lebesgue almost all points, if and only if the
following statement holds:
\renewcommand{\theenumi}{\roman{enumi}}
 \begin{enumerate}
\item \label{countab}
 There exist (up to countable many) ${\cal O}_n \in {\cal O}$
such that $s_n = attrSize ({\cal O}_n) >0$ and $\sum s_n = 1$.

Note that two different trivial reductions  of $\cal O$ are always
mutually independent. Therefore, \ref{countab} is equivalent to the
following:

\item The family ${\cal O}_n$ is a co-chain of trivial reductions
such that $\sum s_n = 1$.
\end{enumerate}
\renewcommand{\theenumi}{\roman{enumi}}
 So the first assertion of Theorem \ref{teoremaCo-cadena} is proved.

The second assertion is trivial from the characterization of physical measures
as those reductions of the space of observable measures that have zero diameter.

To prove the third assertion
observe that each reduction of any co-chain must contain at least
one of the generalized attractors ${\cal O}_n$ because $\sum s_n = 1$, and that
two different reductions of the same co-chain can not contain
 any common reduction, because they must be independent. Then $k$ is
smaller than or equal to the number of independent attractors. On the other
hand, the set of all the independent reductions form itself a
co-chain, so $k$ is greater than or equal to the number of independent attractors. Finally
apply Theorem \ref{teoremaDilworth} to show that
the number of irreducible independent attractors is also  equal to $h$. $\Box$

%\begin{prueba} \em \label{srb} \label{corolarioCo-cadenaphysical}
%{\bf Proof of Theorem \ref{teoremaexistencia0} (b)} It is straightforward consequence
%of Theorems \ref{toeremaminimal0} and \ref{teoremaCardinalidad}. $\Box$
%\end{prueba}

\begin{corollary} \label{corolarioCo-cadenaSRB}
A map $f\colon M \mapsto M$ has \em (up to countable infinitely
many) \em physical measures whose basins of attraction cover
Lebesgue almost all point in $M$, if and only if there exist a \em
(finite or infinite) \em family
$$\{{\cal O}_i ^j, \; i \in I \subset \mathbb{N}, \;  j \in J \subset
\mathbb{N} \}$$ of generalized  attractors ${\cal O}_i^j$ for
$f$ such that \em for all $i \in I$ and $j, k \in J, \; j \neq k $:
$${\cal O}_{i+1}^j \subset  {\cal O}_i^j, \; \;\;
\lim _  {i  } d_i^j \, = \,  0 \; \; \; \lim _{i
 } s_i^{j,k} = 0 \;\;\;\; \mbox{and}
 \;\;\;\lim _{i  } \sum_{j \in J} s_i^j
\, \geq \, 1   $$ where $d_i^j$ and $s_i^j$ denote respectively
the diameter and attracting size of ${\cal O}_i ^j$ and $s_i^{j,k}
$ is the Lebesgue measure of the basin of attraction of ${\cal
O}_i^j \cap {\cal O}_i^k$.

\end{corollary}

{\em Proof: }
If there exist such physical measures $\mu _j$ for $j \in
J \subset N$, simply define the family ${\cal O}_i ^j = \{ \mu
_j\}$ for $i= 1$ and $j \in J$. This family of ergodic attractors
verify all stated conditions.

 To prove the converse statement let us first apply
  Theorem \ref{teoremaCadena} to the chains $\{ {\cal O}_i, i
\in \widehat I \subset I\}$, for each fixed $j \in J$ such that
$\lim _{i } s_i^j
>0$, (while $\lim _i d_ i^j = 0$). Then each of  such chains has a intersection $\{\mu _j \}$
where $\mu _j$ is a  physical measure.

 For $j \neq k$ the basins of
attraction of  $\mu _j$ and $\mu_k$ are Lebesgue almost disjoint,
because its Lebesgue measure is $\lim s_i^{j,k} = 0 $. Thus $\mu
_j \neq \mu _k$. Finally consider the co-chain $\{ \{\mu_j\}, j
\in J\}$ and apply Theorem \ref{teoremaCo-cadena}. $\Box$

\begin{prueba} \em \label{teoremadescomposicionGenArgAttract} {\bf
Proof of Theorem \ref{teoremaDescomposicionGenErgAttr0}:
Decomposition in independent generalized attractors.} We will prove
Theorem \ref{teoremaDescomposicionGenErgAttr0}  in the following
version \ref{teoremaDescomposicionGenAttr1}, that gives an upper
bound to the number of generalized  independent ergodic attractors
in which the space $\cal P$ can be decomposed, and besides states a
sufficient condition for the decomposition be unique. From the
following Theorem \ref{teoremaDescomposicionGenAttr1}, it follows
also the last assertion of Theorem
\ref{teoremaDescomposicionGenErgAttr0}, using that the space ${\cal
O}$ of all observable measures is compact, and thus, for all
$\varepsilon >0$ it can be covered by a finite number of weak$^*$
balls of size $\varepsilon >0$.
\end{prueba}

\begin{theorem} \em {\bf (Reformulation of Theorem \ref{teoremaDescomposicionGenErgAttr0}.)}  \em \label{teoremaDescomposicionGenAttr1}

Any continuous  map $f\colon M \mapsto M$ has  a collection $S$
formed by \em (up to countable many) \em pairwise independent
generalized attractors \em (that are not necessarily irreducible)
\em whose basins of attraction cover Lebesgue almost all points in
$M$.

The supreme $a$ of the  number of such  generalized attractors
verifies ${a} \leq {k} $ \em (where  ${k}$, may be infinite, is the
supreme of  the lengths of the co-chains of reductions in the space
of the observable measures for $f$). \em

If there exists such a  collection $S$  whose generalized attractors
are besides all irreducible, then such $S$ is unique and besides $a
= k = l$, where $l$ is the cardinality of $S$.

\end{theorem}

 {\bf Proof: } To prove the first  and second statements note that,
 by definition of the independence of the reductions, any co-chain $ P$ of reductions
 of the space ${\cal O}$ of the observable measures, verifies
 $\sum s_j \leq 1$, where $s_j$ denotes
  the attracting size of the reduction ${\cal O}_j \in  P$.
  Now take the family ${\cal F}$ of all the co-chains  $S $ such that
   $\sum s_j = 1$.
  (There exists always at
 least one such co-chain: in fact the length-1 co-chain $\{{\cal O}\}$ verifies
 $\sum s_j = s({\cal O }) = 1$, due to the results of Theorem 
 \ref{teoremaminimal}). By construction each $S \in {\cal F}$  verifies the wanted
 conditions.
 As the family ${\cal F}$ is a subfamily of all the co-chains of reductions,
 we obtain $a \leq k$.

  Let us prove now the last assertions of the theorem. If there exists in
  $\cal F$ a co-chain $S= \{\widehat {\cal O}_h: \, h \in H \subset \mathbb{N}\}$
  with cardinality $l = \#H$ and
  whose reductions $\widehat {\cal O}_h$ are all irreducible, to prove that $a = k = l$
  it is enough to
  show than  $l \geq k$.

  Take any
   co-chain $P = \{{\cal O}_ j , \; j \in J \subset \mathbb{N}\} $
   ($P$ is not necessarily in ${\cal F}$).

  It is enough to exhibit an injective application from each $j \in
  J$ to some $h \in H$.

   In fact, let us fix  any $j \in J $   and consider the  basin of attraction
   $C({\cal O}_j)$. By definition of reduction, this basin
  has positive Lebesgue measure $m(C({\cal O}_j))$. But  $S \in {\cal F}$,
   so $\sum _{h \in H} m(C(\widehat
  {\cal O}_h)) = 1$. Then we deduce that  $$0< m (C({\cal O}_j))= \sum
  _{h \in H} m \left( C({\cal O}_j) \cap C(\widehat {\cal O}_h) \right) = \sum
  _{h \in H} m \left( C({\cal O}_j \cap \widehat {\cal O}_h) \right)$$
   Therefore some of the intersections in the sum above at right has positive
  Lebesgue measure. We obtain that for all $j \in J $ there exist some $h = h(j) \in
  H$ such that ${\cal O}_j \cap \widehat {\cal O}_h $ is a
  reduction. But as $\widehat {\cal O}_h $ is irreducible then $ {\cal O}_j \supset \widehat {\cal
  O}_h$.

  To end the proof it remains to show that for $j \neq i \in J$ the
  sets $\widehat {\cal
  O}_{h (i)}$ and  $\widehat {\cal
  O}_{h (j)}$ in $S$ are different reductions. By contradiction, if they were
  the same reduction in $S$, they both
   would be contained in two different
  reductions ${\cal O}_i$ and ${\cal O}_j$ in the chain $P$ and
  therefore these two last reductions would not be independent and
  $P$ would not be a co-chain.

  Let us prove now the unicity of $S$, if there exists one, such
  that
   $S= \{\widehat {\cal O}_h: \, h \in H \subset \mathbb{N}\} \in {\cal F}$
   and
  ${\cal O}_h$ are all irreducible. If there were two such
  collections $S_1$ and $S_2$, then  repeating the
  construction above in this proof with $S = S_1$ and $P = S_2$
  we conclude  any ${\cal O}_j \in S_2$ contains one and
  only one $\widehat {\cal O}_{h(j)} \in S_1$.  But as all $ {\cal O}_j \in S_2$ are
  irreducible we must have $\widehat {\cal O}_{h(j)}= {\cal O}_{j} \in
  S_1$. So $S_2 \subset S_1$. The symmetric relation is obtained
  taking $S= S_2$ and $P = S_1$.
 $\Box$

\section{Appendix.}

\begin{prueba} \em
{\bf Topological attraction in mean of the generalized attractors.}
\end{prueba}

Let $f \in C^0(M)$. Recall Definition \ref{definicionreduccion0} of
Generalized Attractor. Let $(A, {\cal A}) \subset M \times {\cal O}$
be a generalized Attractor, and let $B \subset A$ be its basin of
attraction, which by definition has positive Lebesgue measure $m(B)
>0$. Recall that $A \subset M$ is the minimum compact subset in $M$
that contains the support of all the observable measures $\mu \in
{\cal A}$. In the next statement we will call $A \subset M$ as the
attractor.

\begin{proposition} \em
\label{propositionAtraeEnMedia}

For all  $\varepsilon >0$ there exists $N \geq 1$  and  a subset $C
\subset B$ of the basin of attraction $B$ of the generalized
attractor $A$ such that $m(B\setminus C)< \varepsilon \, m(B)$, and
for all $x \in C$ and for all $n \geq N$, more than
$(1-\varepsilon)100$\% of the iterates of the finite piece
$\{f^j(x)\}_{0 \leq j \leq n-1}$ of the future orbit of $x$, lay  in
the $\varepsilon$-neighborhood of the attractor $A$.

\end{proposition}

{\em Proof:} The attractor $A$ is compact. Call $V \supset A$ to the
$\varepsilon$-neighborhood of $A$ Construct a continuous function
$\phi\colon M \mapsto [0,1]$ such that $\phi(x) = 1$  for all $x \in
A$ and $\phi(x) = 0$ for all $x \not \in V$. By definition of the
generalized attractor, for all $x \in B$ the convergent subsequences
of (\ref{1}) converge to a probability supported in $A$. So $
\frac{1}{n} \sum_{j=0}^{n-1} \phi(f^j(x))$ has all its convergent
subsequences converging to 1. Then:
$$B \subset \bigcup_{N \geq 1} \bigcap_{n \geq N} \{x \in M:
\frac{1}{n} \sum_{j=0}^{n-1} \phi(f^j(x)) > 1 - \varepsilon\}$$
$$m(B) \leq \lim_{N \rightarrow + \infty} m
(\bigcap_{n \geq N} \{x \in M: \frac{1}{n} \sum_{j=0}^{n-1} \phi(f^j(x)) > 1 - \varepsilon\}   \ )$$
Therefore, there exists $N \geq 1$ such that $m(B \setminus C) > 1
-\varepsilon$ where:
$$C =  \bigcap_{n \geq N} \{x \in B: \frac{1}{n} \sum_{j=0}^{n-1} \phi(f^j(x)) > 1 - \varepsilon\}    $$
Call $\chi_V$ to the characteristic function of the  open set $V$.
By construction $0 \leq \chi_V \leq \phi$. Then, for all $n \geq N$,
and for all $x \in C$:
$$\frac{\# \{0 \leq j \leq n-1: f^j(x) \in V\}}{n} = \frac{1}{n} \sum_{j=0}^{n-1} \chi(f^j(x))
 \geq \frac{1}{n} \sum_{j=0}^{n-1} \phi(f^j(x)) > 1 - \varepsilon\; \; \Box$$

\begin{prueba} \em
\label{prueba} {\bf Proof of Theorem \ref{teoremaconvexo}.}
\renewcommand{\theenumi}{\roman{enumi}}
We shall prove the following: \begin{enumerate} \item \label{prime}
If $\mu, \nu \in p\omega(x)$ then for each real number $ 0 \leq
\lambda \leq 1$ there exists a measure $\mu_{\lambda} \in p\omega
(x)$ such that
$$\dist(\mu_{\lambda}, \mu) = \lambda \dist(\mu, \nu)$$
 \item \label{segu} The
set $p\omega(x)$  either has a single element or non-countable
infinitely many. \end{enumerate}
\end{prueba}

\renewcommand{\theenumi}{\arabic{enumi}}
{\em Proof:} First let us deduce \ref{segu} from \ref{prime}:
Suppose that $p\omega (x)$ has at least two different values $\mu$
and $\nu$. It is enough to note that the application $\lambda \in
[0,1] \mapsto \mu_\lambda \in p\omega(x)$ that verifies thesis
\ref{prime}, is injective. Therefore $p\omega(x)$ has non-countable
infinitely many elements.

To prove \ref{prime} consider the sequence (\ref{1}) $ \ \mu_n=
\left\{ \frac1n \; \sum_{j=0}^{n-1} \delta _{f^j(x)}\right\} _{ n
\in \mathbb{N}} $  of time averages. Either it is convergent, or has
at least two convergent subsequences, say $\mu_{m_j} \rightarrow \mu
$ and $\mu _{n_j} \rightarrow \nu $, with  $\mu \neq \nu$. It is
enough to exhibit in the case $\mu \neq \nu$ a convergent
subsequence
 of (\ref{1}) whose limit $\mu_{\lambda } $ verifies the thesis \ref{prime}.

 {\bf  Assertion A: } \em For any given $\varepsilon >0 $ and $K >0$ there exists
 a natural number $h = h(\varepsilon, K)>K$ such that \em
 $$|\dist (\mu_h, \mu) -\lambda \dist (\nu, \mu)| \leq \varepsilon$$

 Let us first prove that Assertion A implies thesis \ref{prime}: Take in
 assertion $A$: $h_0 = 1$ and by induction, for $j \geq 1$ take $h_j$
 given
 $\varepsilon _j = 1/j$ and $K_j = h_{j-1}$. Then we obtain a
 sequence $\mu _{h_j}$, subsequence of (\ref{1}), that verifies $\dist (\mu_{h_j}, \mu)
 \rightarrow \lambda (\dist (\nu, \mu))$. Any convergent
 subsequence of $\mu _{h_j}$ (that do exist  $\cal P$ is compact
 in the weak$^*$ topology) verify \ref{prime}.

  Now, let us prove Assertion A:

As $\mu_{m_j} \rightarrow \mu$ and $\mu_{n_j} \rightarrow \nu $ let
us choose first $m_j$ and then $n_j$ such that
$$m_j >K; \;\; \; \frac{1}{m_j} < \varepsilon /4; \;\; \;  \dist (\mu, \mu_{m_j}) < \varepsilon /4; \; \;\;
  n_j > m_j; \; \;\; \dist (\nu, \mu_{n_j}) < \varepsilon /4$$
  To exhibit the computations let us explicit some metric structure giving
  the weak$^*$ topology of $\cal P$.
  We will use for instance  the following distance: $$\dist (\rho, \delta) =
  \sum_{i=1}^{\infty}\frac {1}{2^i} \; \left| \int g_i \, d \rho
  - \int g_i \, d \delta \; \right|$$ for any $\rho, \delta \in {\cal P}$,
  where $\{g_i\}_{i \in \mathbb{N}}$ is a countable set of functions $g_i \in
  C(M)$ such that $|g_i|\leq 1$,
  dense in the unitary ball of $C(M)$.

 Note from the sequence (\ref{1}) that $|\int g \, d \mu _n - \int g \, d \mu
 _{n+1}| \leq (1/n) ||g||$ for all $g \in C(M)$ and all $n \geq
 1$. Then in particular for $n = m_j + k$, we
 obtain
 \begin{equation} \label{inequ}
  \dist(\mu_{m_j + k}, \mu _{m_j + k + 1})  \leq \frac{1}{m_j} < \varepsilon /4 \;
 \;\; \mbox { for all }
 k \geq 0
 \end{equation}
Now let us choose a natural number $0 \leq k \leq n_j - m_j$ such
that
$$
\left|\dist(\mu_{m_j}, \mu _{m_j + k})- \lambda \dist(\mu_{m_j}, \mu
_{n_j}) \right|< \varepsilon /4 \;\;\mbox{ for the given } \lambda
\in [0,1]$$ Such $k$ does exist because inequality (\ref{inequ}) is
verified for all $k \geq 0$ and besides:

$\bullet$ If $ k= 0 \mbox{ then } \dist(\mu_{m_j}, \mu _{m_j + k})=
0$ and

 $\bullet$ if $ k= n_j -m_j \mbox{ then }  \dist(\mu_{m_j}, \mu _{m_j + k})=
\dist(\mu_{m_j}, \mu _{n_j })$

Now renaming $h = m_j + k$, joining all the inequalities above, and
applying the triangular property, we deduce:
$$\left| \dist(\mu_{h}, \mu) - \lambda \dist (\nu, \mu)\right| \leq
 \left| \dist(\mu_{h}, \mu_{m_j}) - \lambda \dist (\mu_{m_j},
\mu_{n_j})\right| + $$ $$+ \left| \dist(\mu_{h}, \mu) -
\dist(\mu_{h}, \mu_{m_j})\right| + $$ $$ + \lambda \left| \dist
(\mu_{m_j}, \mu_{n_j}) - \dist (\mu_{m_j}, \nu)\right| + \lambda
\left| \dist (\mu_{m_j}, \nu) - \dist (\mu, \nu)\right| <
\varepsilon \;\;\;\;\; \Box
$$

{\bf Acknowledgements.}

We thank the referee for his   remarks and suggestions to the first version of this paper.

\vspace{.4cm}

%$\clubsuit$ Sacar de la lista que sigue las referencias bibliogr\'{a}ficas que no son usadas en el trabajo.

\end{document}